\documentclass[11pt,letterpaper]{amsart}
\pdfoutput=1
\usepackage[T1]{fontenc}
\usepackage{lmodern,amssymb,mathrsfs}
\usepackage{mathtools}
\usepackage{color}
\usepackage[hmargin=2.5cm,vmargin=2cm]{geometry}
\usepackage[normalem]{ulem}
\usepackage{tcolorbox}
\usepackage{mathabx}
\usepackage{hyperref}
\usepackage{float}
\hypersetup{
colorlinks=true,
linkcolor={blue},
citecolor={red!50!black}
}

\newtheorem{theorem}{Theorem}[section]
\newtheorem{maintheorem}{Theorem}
\newtheorem*{maintheoremrepeat}{Theorem~\ref*{thm:main}}

\newtheorem{proposition}[theorem]{Proposition}
\newtheorem{lemma}[theorem]{Lemma}
\newtheorem{corollary}[theorem]{Corollary}
\theoremstyle{definition}
\newtheorem{definition}[theorem]{Definition}
\newtheorem{example}[theorem]{Example}
\theoremstyle{remark}

\numberwithin{equation}{section}

\newcommand{\C}{\mathbb C}
\newcommand{\Nz}{\mathbb N}
\newcommand{\B}{\mathbb B}
\newcommand{\cF}{\mathcal F}

\newcommand{\cS}{\mathcal S}
\newcommand{\Ann}{\operatorname{Ann}}

\title[A proof of the Arveson--Douglas conjecture]{A proof of the Arveson--Douglas conjecture}

\author[C. Bruce]{Chris Bruce}\thanks{Chris Bruce was supported by the Engineering and Physical Sciences Research Council, grant no. UKRI3684.}
\address[C. Bruce]{School of Mathematics, Statistics and Physics, Herschel Building, Newcastle University, Newcastle upon Tyne, NE1 7RU, UK}
\email{chris.bruce@newcastle.ac.uk}

\author[A. Dor-On]{Adam Dor-On} \thanks{Adam Dor-On was partially
supported by an ISF grant no. 2499/25, by an NSF-BSF grant no. 2350543 / 2023695 (respectively), and a DFG Middle-Eastern collaboration project no. 529300231.}
\address[A. Dor-On]{Department of Mathematics, University of Haifa, Mount Carmel, Haifa 3103301, Israel}
\email{adoron.math@gmail.com}

\date{\today}



\begin{document}

\begin{abstract}
We provide a proof of the Arveson--Douglas $p$-essential normality conjecture. More precisely, for a homogeneous ideal in the ring of complex polynomials in $d$ commuting variables with multiplicity $r$, we show that when the homogeneous variety of the annihilator of the quotient module has dimension $\delta$, the commutators $[S_i^*, S_j]$ belong to the Schatten class $\mathcal S_p$ for every $p>\delta$. 
\end{abstract}

\maketitle

\vspace{-18pt}

\section{Introduction}\label{sec:target}

The coordinate shifts on the Drury--Arveson space commute, but they do not commute with each other's adjoints. Their failure to do so is measured by compact operators, and more precisely by operators in Schatten classes. Arveson's essential normality conjecture asks whether this property survives when homogeneous polynomial relations are imposed on the shifts. Douglas proposed that the dimension of the variety of the homogeneous ideal determines the range of Schatten exponents. In this paper, we provide an accessible proofs of the finite-multiplicity form of this conjecture.
 
For a compact operator $A$ we write $s_1(A)\ge s_2(A)\ge\cdots$ for its singular values, namely the eigenvalues of $|A|=(A^*A)^{1/2}$, listed with multiplicity and continued by zeros when $A$ has finite rank. For $p>0$, the Schatten class $\cS_p$ consists of the compact operators $A$ for which $\sum_{k\ge1}s_k(A)^p<\infty$.

\begin{definition}\label{def:essential-normality}
Let $T=(T_1,\ldots,T_d)$ be a commuting tuple of bounded operators on a Hilbert space. We say that $T$ is \emph{essentially normal} if $[T_i^*,T_j]$ is compact for every $1\le i,j\le d$. For $p>0$ we say that $T$ is \emph{$p$-essentially normal} if $[T_i^*,T_j]\in\cS_p$ for every $1\le i,j\le d$.
\end{definition}

Let $\B_d$ denote the open unit ball of $\C^d$. The Drury--Arveson space
$H_d^2$ is the reproducing-kernel Hilbert space on $\B_d$ with kernel
$k(z,w)=(1-\langle z,w\rangle)^{-1}$, where
$\langle z,w\rangle=\sum_j z_j\overline{w_j}$. Thus,
$f(w)=\langle f,k(\cdot,w)\rangle_{H_d^2}$. We write $M_{z_j}$ for
multiplication by $z_j$ on $H_d^2$. For $d,r\ge1$, put
$R_d=\C[z_1,\ldots,z_d]$, and consider $R_d^r$ as a module over $R_d$. A submodule
$J\subseteq R_d^r$ is homogeneous if it contains every homogeneous component of each of its elements. We write $J_n$ for its degree-$n$ part.

Quotienting $H_d^2\otimes\C^r$ by the closure of $J$ imposes the relations in $J$ on coordinate multiplication. To compute the adjoints and commutators of the induced operators, we realize each coset by its unique representative orthogonal to that closure. In what follows, we denote by $P_E$ the orthogonal projection onto a closed subspace
$E$.

\begin{definition}\label{def:quotient-shift}
Let $J\subseteq R_d^r$ be homogeneous and let $\overline J$ denote its
closure in $H_d^2\otimes\C^r$. The \emph{quotient module} associated
with $J$ is represented on
\[
  \mathcal Q_J=(H_d^2\otimes\C^r)\ominus\overline J .
\]
For $1\le j\le d$, put $T_{z_j}=M_{z_j}\otimes1_{\C^r}$. The
\emph{quotient shift} on $\mathcal Q_J$ is the tuple
$S=(S_1,\ldots,S_d)$ given by $S_j=P_{\mathcal Q_J}T_{z_j}|_{\mathcal Q_J}$. We write $\mathcal Q_{J,n}$ for the space of homogeneous vectors
of degree $n$ in $\mathcal Q_J$.
\end{definition}

Since each $T_{z_j}$ preserves $\overline J$, the compression in
Definition~\ref{def:quotient-shift} represents multiplication by $z_j$
on cosets, so the quotient shifts commute. For the algebraic quotient $M=R_d^r/J$, its annihilator and its complex variety are given by
\[
 I=\Ann M=\{p\in R_d:pR_d^r\subseteq J\},\ \ \ 
 V(I)=\{z\in\C^d:p(z)=0\text{ for all }p\in I\}.
\]
When $r=1$ one has $I=J$. For $M\ne0$ we write $\delta=\dim_\C V(I)$ for the affine dimension of the variety.
The summability threshold will use this dimension, with no radicality assumption on $I$.

Arveson asked whether every graded pure $d$-contraction of finite rank
is $p$-essentially normal for $p>d$ \cite[Problem~2]{arveson2002}. He subsequently formulated the conjecture for homogeneous quotients of the finite-multiplicity $d$-shift \cite[Conjecture~1]{arveson2003problems}. For scalar homogeneous quotients, Douglas proposed the dimension of the variety as the threshold for summability \cite[discussion following Theorem~5.3]{douglas2006}. The finite-multiplicity formulation asks for $p>\delta$ \cite[p.~276 and \S4]{shalitdivision}. Our main result gives a bound on each homogeneous degree from which we get a decay for the singular values.

\begin{maintheorem}\label{thm:main}
For every homogeneous submodule $J\subseteq R_d^r$ of infinite
codimension and every
$1\leq i,j\leq d$, there is a constant $C_{ij}$, independent of $n$,
such that for every $n\in \mathbb{N}$,
\begin{equation}\label{eq:main-bound}
\bigl\|[S_i^*,S_j]|_{\mathcal Q_{J,n}}\bigr\|
\leq\frac{C_{ij}}{n+1}.
\end{equation}
Writing $I=\Ann(R_d^r/J)$ and $\delta=\dim_\C V(I)$, the singular values satisfy $s_k([S_i^*,S_j])=O(k^{-1/\delta})$, so it follows that $[S_i^*,S_j]\in\cS_p$ for every $p>\delta$.
\end{maintheorem}

Theorem~\ref{thm:main} also answers Arveson's original question \cite[Problem~2]{arveson2002} about graded pure row contractions of finite rank, and the reduction to the quotient formulation is proven in Corollary \ref{cor:finite-rank-row-cont}.

The strict inequality $p>\delta$ in Theorem~\ref{thm:main} is necessary in general for every positive dimension. The case $J=0$ already illustrates the role of dimension in this problem. Indeed, Arveson proved that for $d\geq2$ the $d$-shift $(M_{z_1},\ldots,M_{z_d})$ on $H_d^2$ is $p$-essentially normal exactly in the range $p>d$ \cite[Proposition~5.3]{arveson1998}.

For quotients of $H_d^2\otimes\C^r$ by submodules generated by monomials,
Arveson proved $p$-essential normality for $p>d$
\cite[Corollary~2.2]{arveson2005}, and Douglas obtained the improved range $p>\dim_\C V(I)$ \cite[Theorem~5.3]{douglas2006}. Guo and K.~Wang proved that homogeneous principal submodules of $H_d^2$ are $p$-essentially normal for $p>d$ in every dimension \cite[Theorem~2.2]{guowang}. For scalar homogeneous ideals they obtained the dimension-dependent range in two and three variables, and in any number of variables for an ideal $gL$ with a common factor $g$ of positive degree where $\dim_\C V(L)\leq1$ \cite[Propositions~4.2--4.3]{guowang}.
Their finite-multiplicity results give $p>2$ in two variables and $p>1$ in any number of variables when the variety is one-dimensional \cite[Theorem~3.1 and Proposition~4.1]{guowang}. Bauer later gave a Fock-space proof of their theorem on principal submodules \cite[Theorem~2.7 and Example~2.11]{bauer2026}.

The geometric form of the conjecture \cite[Conjecture~1.1]{kennedyshalit2012} predicts that, for a homogeneous
variety $V\subseteq\C^d$, the quotient of $H_d^2$ by the functions
vanishing on $V\cap\B_d$ is $p$-essentially normal for
$p>\dim_\C V$. This is the Arveson--Douglas conjecture for radical homogeneous ideals $J\subseteq R_d$. The identification of the quotient with the closed span of the reproducing kernels on the variety is precisely the radical ideal criterion of Davidson, Ramsey, and Shalit
\cite[Lemma~5.5]{davidsonramseyshalit2011}. Using Hartz's closed-sum
theorem \cite[Corollary~5.8]{hartz2012} and \cite[Theorem~3.6]{kennedyshalit2012}, Kennedy and Shalit proved this geometric form of the conjecture for finite unions of linear subspaces \cite[Theorem~4.2]{kennedyshalit2012}. In a breakthrough paper, Engli\v{s} and Eschmeier proved the geometric conjecture for homogeneous varieties smooth away from the origin, obtaining the singular-value estimate in Theorem~\ref{thm:main} as well \cite[Main Theorem, p.~627, and
Theorem~5]{engliseschmeier}. Following this, Ding obtained the range $p>\dim_\C(V\cap\B_d)$ for analytic varieties $V$ defined near the closed ball that are smooth along the sphere and meet it
transversally \cite[Corollary~1.5(1) and Lemma~3.4(4)]{ding2026}.
In the case of a pure one-dimensional analytic variety with these
boundary hypotheses, Xia proved that the commutators are trace class
\cite[Theorem~1.1]{xia2024}.

Other approaches to the homogeneous conjecture include Kennedy's
decomposability criterion \cite{kennedy2015} and the use of Veronese
embeddings in the geometric setting \cite{aricige2026}. For surveys on the Arveson--Douglas conjecture, we recommend
\cite{guoyiwangsurvey}, \cite[\S10.3.5]{hartz2023}, and
\cite{hartzshalit2026}.

Several consequences now follow from prior works, and we mention but a few in this version of the paper. For $p>d$, Arveson's restriction and compression theorem gives $p$-essential normality of the coordinate shifts restricted to $\overline J$ \cite[Theorem~4.3(ii)--(iv)]{arveson2005}. For an infinite-dimensional quotient, the images of the quotient shifts modulo compact operators generate $C(V(I)\cap\partial\B_d)$. For the scalar case, this follows for example from \cite[Theorem~5.1(4)]{guowang}. In finite multiplicity, the identification follows by combining the interior fiber description in \cite[Corollary~4]{didaseschmeierhartzscherer2024} with the boundary spectrum criterion in \cite[Corollary~3.3]{wernet2014}. The resulting homomorphism from $C(V(I)\cap\partial\B_d)$ into the Calkin algebra then determines an odd $K$-homology class on $V(I)\cap\partial\B_d$ \cite{arveson2005, bdf1977}.

Two approaches to the conjecture are particularly relevant to our proof.
Shalit showed that $(S_1,\ldots,S_d)$ is $p$-essentially normal for
$p>\dim_\C V(I)$ whenever $J$ admits stable polynomial division in the
Drury--Arveson norm \cite[Theorem~4.1]{shalitdivision}. His estimate
\cite[equation (4.1)]{shalitdivision} bounds the component of the shift from the quotient to the submodule on degree $n$ by $O((n+1)^{-1/2})$.
The compression formula thus gives the degree estimate of
Theorem~\ref{thm:main} in that setting. Eschmeier related essential normality to an estimate involving the number operator, which acts as multiplication by $n$ in degree $n$ \cite{eschmeier2011}. For homogeneous generators of the same degree, the corresponding derivative inequality is recorded in \cite[equation (6.1)]{shalitdivision}. The degree estimate in Theorem~\ref{thm:main} itself is known in several cases. For a homogeneous principal ideal, Guo and K.~Wang factored the projection commutator through $(N+1)^{-1/2}$, where $N$ is the number operator \cite[Corollary~2.5]{guowang}. Their proofs in two variables and for one-dimensional varieties also give degree estimates in finite multiplicity \cite[Theorem~3.1 and Proposition~4.1]{guowang}. Wernet formulated the corresponding factorization for varieties of dimension at most one in any number of variables, and for arbitrary varieties in two variables. In the scalar case he also treated ideals for which the common zero set of the generators after division by their greatest common divisor has dimension at most one at the origin. This includes all homogeneous ideals in at most three variables \cite[Theorems~2.21--2.22, 2.25 and Corollary~2.26]{wernet2014}. In each case the factorization gives the degree estimate by the compression formula, and counting dimensions of the homogeneous pieces gives the stated decay of the singular values. We repeat these two steps in the preliminaries Subsection~\ref{subsec:graded-criterion} for posterity.

To prove Theorem~\ref{thm:main}, we use the fact that on homogeneous
polynomials of a fixed degree the Drury--Arveson and Fock norms differ
only by a scalar. Consider first a homogeneous ideal $I$, and take
orthogonal complements within each degree in Fock space. Multiplication
by $z_j$ on degree $n$ has norm $\sqrt{n+1}$. We shall prove that its
component from $I_n^\perp$ to $I_{n+1}$ has norm bounded independently
of $n$. To see why this suffices, observe that in Drury--Arveson space
the same component acquires a factor $1/\sqrt{n+1}$. In the compression
formula, two such components occur as a product, one of them adjointed.
Their contribution to the quotient commutator in degree $n$ is therefore
$O((n+1)^{-1})$. The product rule for differentiation gives the same
bound for the remaining term. Subsection~\ref{subsec:graded-criterion}
makes this calculation explicit and uses the Hilbert polynomial to
count the singular values contributed by the different degrees.
The proof of the uniform Fock bound is by compactness and does not
compute the acutal bound. The precise dependence on $n$ comes from the
change of norms, so no rate of convergence in the compactness
argument is needed.

Suppose that the proposed uniform bound in Fock space fails. Then, along a sequence of degrees $n\to\infty$ there are unit singular vectors
$u_n\perp I_n$ and $v_n\in I_{n+1}$ with singular values
$\lambda_n\to\infty$. Since $\partial_j=\partial/\partial z_j$ is the
Fock adjoint of multiplication by $z_j$ on polynomials, we have
\[
 z_j u_n-\lambda_n v_n\perp I_{n+1},\qquad
 \partial_jv_n-\lambda_n u_n\in I_n .
\]
Our aim is to pass these two relations to a limit that retains both
ideal membership and the differential equations for orthogonality. The
vectors themselves tend to zero locally uniformly, since the Fock norm
gives $|u_n(z)|\leq |z|^n/\sqrt{n!}$ and the corresponding bound in
degree $n+1$ holds for $v_n$. Instead, we translate the sum $u_n+v_n$
to a point where its Gaussian-weighted absolute value is largest, and
normalize the translate to have value one at the origin. The two terms
have consecutive degrees. Rotating $z$ by a common scalar of modulus
one therefore allows their values to have the same argument. This
observation bounds both translated terms with the same normalization.

Translation moves the ideal as well. In Section~\ref{sec:affine-limits}
we homogenize the translated ideals and pass to a limit in each degree.
Let $L$ denote the affine ideal obtained by setting the homogenizing
variable equal to one in the limit. Two conditions on entire functions pass to this limit. The first is local uniform approximation by ideal polynomials. Here, the approximating polynomials can have arbitrarily large degrees. We express them in terms of finitely many ideal polynomials using coefficient estimates independent of the degree.
The second condition is the system of differential equations
$p^*(\partial)g=0$ for $p\in L$, where for
$p(z)=\sum_\alpha p_\alpha z^\alpha$ and
$\partial^\alpha=\prod_j\partial_j^{\alpha_j}$ and we define
$p^*(\partial)=\sum_\alpha\overline{p_\alpha}\partial^\alpha$.
On Fock space these equations express orthogonality to $L$.
After canceling the translation-center terms, dividing the two
translated singular-vector relations by $\lambda_n$ shows that a
limit of the normalized sums satisfies both conditions on entire functions.

The uniqueness theorem in Section~\ref{sec:affine-uniqueness}
shows that such a function vanishes for every affine ideal when its
Gaussian-weighted absolute value has polynomial growth. If it decays faster than every
inverse power, we write the function as a sum of polynomials from the ideal
with entire coefficients. H\"ormander's division theorem allows us to
choose these coefficients with a weighted square-integrability estimate.
In the Fock pairing of the function with itself, we can therefore
truncate the coefficients. Each resulting pairing vanishes by the
differential equations. If this decay condition fails, we translate once more, along a sequence tending to infinity, and obtain
a function with value one at the origin. The new limiting ideal is
mapped into itself by a nonzero translation. Restricting the
function to the perpendicular hyperplane preserves its value at the
origin and both ideal membership conditions. The normalization may increase the polynomial growth exponent, which is why uniqueness is proved for arbitrary such exponent at once. Induction on the number of variables therefore proves the required vanishing. In Section~\ref{sec:singular-vectors} we apply this to the normalized sum
and use Shalit's realization of vector-valued polynomials as scalar
polynomials linear in auxiliary variables \cite[\S5]{shalitdivision}
to obtain the Fock bound in every finite multiplicity.

\smallskip
\noindent \textbf{Acknowledgments.} The second-named author is grateful to Camila Sehnem, for discussions about the Arveson--Douglas conjecture held at the Research Institute for Mathematical Sciences, Kyoto University. Both authors thank Kevin Aguyar Brix for stimulating conversations related to the Arveson--Douglas conjecture that took place at a 2025 Research in Groups meeting at the ICMS in Edinburgh, as well as thank Ian Thompson for discussions held during the workshop ``Around the Arveson--Douglas Conjecture'' at the University of Potsdam. 

\smallskip
\noindent \textbf{AI disclosure.} GPT-5.6 Sol and GPT-6.0 Astra from OpenAI, as well as Claude Opus 5.0 and Claude Opus 5.5 from Anthropic were used to search for strategies involving the use of group symmetry, as well as boundary limiting techniques from harmonic analysis, in developing our initial proof. GPT-6.0 Astra, GPT-6.1 Sol and GPT-6.2 Astra from OpenAI, as well as Claude Opus 5.5 and Claude Fable 5.2 from Anthropic were then used to suggest simplifications and improvements to the the proofs, produce several drafts of the paper, and to proofread the final manuscript. The authors subsequently checked and rewrote the arguments to make them more accessible and take full responsibility for the validity of all claims and proofs in the manuscript.

Chris Bruce is supported by a ChatGPT for Academic Researchers plan and a Claude Team plan for Scientists.

\section{Preliminaries}\label{sec:preliminaries}

We first reduce the commutator estimates to a uniform bound on the
component of coordinate multiplication that sends a quotient graded component into the submodule in Fock space. The following subsections prepare for the proof of that bound. They express ideal membership and Fock orthogonality as conditions on entire functions, describe their
behavior under translation, and introduce the algebraic compactness
and coefficient norms used to pass to limiting ideals.

We write $\Nz=\{0,1,2,\ldots\}$ and take Hilbert space inner products to be linear in the first variable. We retain the notation for projections,
commutators, and Schatten classes from the introduction.

\subsection{Drury--Arveson and Fock spaces}\label{subsec:multivariable-operators}

The two Hilbert spaces used in the proof have the same orthogonal
monomials, but different normalizations. We record how this change
of scale affects graded components of the quotient and how it maps between successive degrees.

Let $R_{d,n}$ be the space of homogeneous polynomials of degree $n$ in
$R_d$. For a homogeneous submodule $J\subseteq R_d^r$, put
$J_n=J\cap R_{d,n}^r$, so that algebraically we have $J=\bigoplus_{n\ge0}J_n$. For $\alpha\in\Nz^d$ write $|\alpha|=\sum_j\alpha_j$, as well as
$\alpha!=\prod_j\alpha_j!$ and $z^\alpha=\prod_jz_j^{\alpha_j}$. The added level of symmetry needed for the proof arises from considering weighted translations on Fock space of entire functions $f$ on $\C^d$ which have finite norm
\[
 \|f\|_{\cF_d}^2=\int_{\C^d}|f(z)|^2e^{-|z|^2}\,d\mu(z),
 \qquad d\mu(z)=\pi^{-d}\,dz,
\]
where $dz$ is Lebesgue measure. For $\C^r$-valued functions, as well as the absolute value of a vector $z\in \mathbb{C}^d$, we mean the Euclidean norms in the integral formula above. We write $\|\cdot\|_{\cF}$ when the number of variables and the multiplicity are understood.
The monomials are orthogonal
in both Drury--Arveson space, and Fock space, with
\begin{equation}\label{eq:two-monomial-norms}
 \|z^\alpha\|_{H_d^2}^2=\frac{\alpha!}{|\alpha|!},
 \ \ \  \|z^\alpha\|_{\cF_d}^2=\alpha!.
\end{equation}
For the first identity, expand the Drury--Arveson kernel
\cite{drury1978,arveson1998}, whose degree-$n$ terms are
$\sum_{|\alpha|=n}(n!/\alpha!)z^\alpha\overline{w^\alpha}$.
For the second, integrate the Taylor series over each of the coordinate tori. Angular integration eliminates mixed monomials, and the
remaining nonnegative terms can be integrated in the radii. The
identity $2\int_0^\infty r^{2a+1}e^{-r^2}\,dr=a!$ then gives
\[
 \Big\|\sum_\alpha c_\alpha z^\alpha\Big\|_{\cF_d\otimes\C^r}^2
 =\sum_\alpha\|c_\alpha\|_{\C^r}^2\alpha!.
\]
Thus, Taylor polynomial truncations converge in Fock norm, and on $R_{d,n}^r$ the Fock inner product is $n!$ times the Drury--Arveson inner product (see \cite[Example~2.3]{bauer2026}).

The same monomial formula explains why we will need weighted translations
when taking limits. For a homogeneous polynomial $f$ of degree
$n$, Cauchy--Schwarz and the multinomial identity give
\begin{equation}\label{eq:homogeneous-evaluation}
 |f(z)|^2\le\|f\|_{\cF_d}^2
       \sum_{|\alpha|=n}\frac{|z^\alpha|^2}{\alpha!}
 =\|f\|_{\cF_d}^2\frac{|z|^{2n}}{n!}.
\end{equation}
Consequently, by Stirling's approximation, a unit vectors in degrees tending to infinity converge to zero on every fixed compact set.

Write $\cF_{d,n}$ for $R_{d,n}$ with the Fock norm. Since $J_n$
is finite dimensional, projection to degree $n$ takes the Hilbert
space closure of $J$ onto $J_n$. The two inner products differ by a
scalar on this graded component, so they give the same polynomial quotient space
\[
 K_n=(\cF_{d,n}\otimes\C^r)\ominus J_n.
\]
As a space of polynomials, $K_n$ is also $\mathcal Q_{J,n}$, but this
identification does not preserve the two Hilbert space norms. In the
Drury--Arveson norm, $\mathcal Q_J=\bigoplus_{n\ge0}K_n$
orthogonally, and sending $f\in K_n$ to $f+J_n$ identifies $K_n$ linearly with $R_{d,n}^r/J_n$.

All Fock adjoints below are taken between finite polynomial graded components. The monomial identities show that multiplication by $z_j$ from degree $n$ to degree $n+1$ has Fock adjoint $\partial_j$ and norm
$\sqrt{n+1}$. In Drury--Arveson norm its adjoint on degree $n+1$
is $\partial_j/(n+1)$, and the backward shift vanishes on constants.
These formulas act componentwise for multiplicity $r$. Moreover,
coordinate multiplication is contractive in Drury--Arveson space,
since its weights on normalized monomials are
$\sqrt{(\alpha_j+1)/(|\alpha|+1)}\le1$. Its quotient shifts $S_j$
are therefore contractions as well.

\subsection{A criterion for the graded commutators}\label{subsec:graded-criterion}

For $f\in K_n$, the product $z_jf$ belongs to the full homogeneous space of degree $n+1$, which decomposes orthogonally as $K_{n+1}\oplus J_{n+1}$. Its component in $K_{n+1}$ is $S_jf$. We shall estimate its component in $J_{n+1}$. For a fixed homogeneous $J\subseteq R_d^r$, define
\[
 B_{j,n}\colon K_n\longrightarrow J_{n+1},\qquad
 B_{j,n}f=P_{J_{n+1}}(z_jf),\qquad
 b_j=\sup_{n\ge0}\|B_{j,n}\|_{\cF}.
\]
Here, the block norm uses the Fock norms on both its domain and its
range, and $b_j$ could apriori be infinite. Its Fock adjoint is
$B_{j,n}^*v=P_{K_n}\partial_jv$. The same polynomial map, regarded
between Drury--Arveson graded components, has norm
$\|B_{j,n}\|_{\cF}/\sqrt{n+1}$, because the domain and range Fock
norms squared have scale factors $n!$ and $(n+1)!$, respectively.
We keep the notation $B_{j,n}^*$ for the Fock adjoint and $S_i^*$ for
the Drury--Arveson adjoint. In particular, the factor $n+1$ in the
next identity accounts for the change of norm between successive
degrees. With these conventions, Arveson's
compression identity \cite[(4.2)]{arvesonquotients} has the following
form on each graded component.

\begin{lemma}\label{lem:graded-compression}
 The commutator $[S_i^*,S_j]$ preserves each $K_n$, and
 \begin{equation}\label{eq:compression}
 (n+1)[S_i^*,S_j]|_{K_n}
 =\delta_{ij}1_{K_n}-S_jS_i^*|_{K_n}-B_{i,n}^*B_{j,n}.
 \end{equation}
 If $b_i,b_j<\infty$, then
 \begin{equation}\label{eq:block-bound}
 \bigl\|[S_i^*,S_j]|_{K_n}\bigr\|
 \le\frac{1+b_i b_j}{n+1}\qquad(n\ge0).
 \end{equation}
\end{lemma}
\begin{proof}
 The inclusion $z_iJ_n\subseteq J_{n+1}$ and adjunction imply that
 the Drury--Arveson backward shift takes $K_{n+1}$ into $K_n$.
 Hence, $S_i^*=\partial_i/(n+1)$ on $K_{n+1}$ and $S_i^*=0$ on $K_0$.
 For $f\in K_n$ we have $S_jf=z_jf-B_{j,n}f$, so
 \[
 (n+1)S_i^*S_jf
 =P_{K_n}\partial_i(z_jf)-B_{i,n}^*B_{j,n}f.
 \]
 The product rule turns the first term into  $\delta_{ij}f+P_{K_n}(z_j\partial_if)$. When $n\ge1$, the second summand of the latter expression is precisely $nS_jS_i^*f$. When $n=0$, both that summand and $S_jS_i^*f$ vanish. Subtracting $(n+1)S_jS_i^*f$ then proves the identity in equation  \eqref{eq:compression}, including for degree zero.

 If $i\ne j$, contractivity gives $\|S_jS_i^*\|\le1$. For $i=j$,
 the inequalities $0\le S_iS_i^*\le1$ give
 $\|1-S_iS_i^*\|\le1$. Thus, the first two terms on the right of equation
 \eqref{eq:compression} have norm at most one. The remaining product $B_{i,n}^*B_{j,n}$ has Fock norm at most $b_i b_j$, and its domain and co-domain are both $K_n$. On $K_n$, the Fock norm is $\sqrt{n!}$ times the Drury--Arveson norm. Since $B_{i,n}^*B_{j,n}$ maps $K_n$ into itself, this factor cancels in the ratio defining its operator norm. Thus, we get that the Fock operator norm of $B_{i,n}^*B_{j,n}$ as an operator on $K_n$ is equal to its Drury--Arveson operator norm, and we have proven the inequality in equation \eqref{eq:block-bound}.
\end{proof}

This is the sufficient direction of Eschmeier's criterion
\cite[Theorem~2.4]{eschmeier2011}, and in finite multiplicity it is
\cite[Theorems~2.19 and~2.21(ii)]{wernet2014}. To count singular
values, we combine the block bound with the dimensions of the
original quotient graded components. For $M=R_d^r/J$, put $M_n=R_{d,n}^r/J_n$. Its Hilbert function is given by $n\mapsto \dim_\C M_n=\dim K_n$. If $M$ is finite dimensional, only finitely
many $K_n$ are nonzero, so all commutators have finite rank and
belong to $\cS_p$ for $p>0$.

Suppose henceforth that $M$ is infinite dimensional. The affine dimension $\delta=\dim_\C V(\Ann M)$ equals the Krull
dimension of $R_d/\Ann M$, that is, the supremum of the lengths of
strict chains of prime ideals in this ring \cite[p.~225]{eisenbud}. This dimension is positive, because a zero-dimensional affine $\C$-algebra is finite dimensional, as would then be its finitely generated module $M$. The Hilbert function of $M$ eventually agrees with a polynomial of degree $\delta-1$ \cite[Lemma~3.3]{guowang}. Enlarging a constant to cover finitely many earlier graded component gives some $c_M\ge1$ such that for all $n\ge0$,
\begin{equation}\label{eq:hilbert-growth}
 \dim K_n\le c_M(n+1)^{\delta-1}.
\end{equation}

The block estimate and the Hilbert function have different roles:
the former bounds the size of each singular value in a given degree,
and the latter bounds how many that degree can contribute.

Shalit uses the Hilbert function in this way to deduce Schatten
summability from a bound on the homogeneous blocks \cite[proof of Theorem~4.1 and Remark~4.5]{shalitdivision}. We record the count for individual singular values, since it also gives the decay asserted in Theorem~\ref{thm:main}.

\begin{proposition}\label{prop:graded}
 Let $J\subseteq R_d^r$ be a homogeneous submodule of infinite
 codimension. If $b_j<\infty$ for every $j$, then
 Theorem~\ref{thm:main} holds for $J$ with $C_{ij}=1+b_i b_j$.
\end{proposition}
\begin{proof}
 On the Drury--Arveson orthogonal sum $\bigoplus_nK_n$, write
 $A=[S_i^*,S_j]=\bigoplus_nA_n$. Lemma~\ref{lem:graded-compression}
 gives $\|A_n\|\le C_{ij}/(n+1)$. Truncation to the first finitely
 many grades therefore approximates $A$ in operator norm by
 finite-rank operators, proving compactness of $A$. Also
 $|A|=\bigoplus_n|A_n|$, so the positive singular values of $A$ are
 those of its finite-dimensional blocks, with multiplicity.

Put $C=C_{ij}=1+b_i b_j$. Fix $k\ge1$ and suppose that $s=s_k(A)>0$. There are at least $k$ singular values of $A$ greater than or equal to $s$. If $A_n$ contributes any such singular value, then $s\le\|A_n\|\le C/(n+1)$, so $n+1\le C/s$. Putting $m=\lfloor C/s\rfloor$, all these singular values therefore come from $A_0,\ldots,A_{m-1}$. Since each of these blocks contributes at most its dimension, equation \eqref{eq:hilbert-growth} gives
\[
 \begin{aligned}
 k
 \le \sum_{n=0}^{m-1}\dim K_n
 \le c_M\sum_{\ell=1}^{m}\ell^{\delta-1}
 \le c_Mm^\delta
 \le c_M(C/s)^\delta.
 \end{aligned}
\]
Here, $\delta\ge1$, so each of the $m$ summands
$\ell^{\delta-1}$ is at most $m^{\delta-1}$. Thus, we get that $s_k(A)\le D \cdot k^{-1/\delta}$ for $D = C \cdot c_M^{\frac{1}{\delta}}$, which also holds when $s_k(A)=0$. Consequently, $\sum_{k\ge1}s_k(A)^p<\infty$ for every $p>\delta$,
by first comparison test with the convergent series $\sum_{k\ge1}k^{-p/\delta}$.
\end{proof}

\subsection{Affine ideals and Fock orthogonality}\label{subsec:affine-spaces}

To bound $b_j$ for a homogeneous ideal $I$, we shall translate
polynomials in $I$ and polynomials orthogonal to $I$ and take locally
uniform limits in Section~\ref{sec:singular-vectors}. We formulate the
corresponding conditions on entire functions using compact-open
approximation and constant-coefficient differential equations. These
conditions still make sense when a limiting function has possibly infinite Fock norm. Translation need not preserve homogeneity, so we formulate both conditions for arbitrary polynomial ideals. Weighted translation will carry both conditions to the translated ideal while translating the Gaussian-weighted absolute value.

We call an ideal of $R_d$ an \emph{affine ideal} to emphasize that no homogeneity assumption is made on it. Equip the space
$\mathcal O(\C^d)$ of entire functions with the topology of uniform
convergence on compact sets, and for $\beta\in\Nz^d$ write $\partial^\beta=\prod_j\partial_j^{\beta_j}$.

\begin{definition}\label{def:affine-function-spaces}
 Let $I\subseteq R_d$ be an affine ideal. Write $\mathcal A(I)$
 for its \emph{compact-open closure} in $\mathcal O(\C^d)$. For a polynomial
 $p=\sum_\alpha p_\alpha z^\alpha$ let
 $p^*(\partial)=\sum_\alpha\overline{p_\alpha}\,\partial^\alpha$,
 and define
 \[
 \mathcal N(I)=\{g\in\mathcal O(\C^d):p^*(\partial)g=0
                              \text{ for every }p\in I\}.
 \]
\end{definition}

Thus, $\mathcal A(I)$ yields entire functions which can be approximated by polynomials in the ideal, whereas $\mathcal N(I)$ is described via a system of differential equations. Both are complex vector spaces. Their definitions retain the ideal $I$ itself, including its multiplicities, rather than merely its zero set.

The conjugates in $p^*(\partial)$ are dictated by Fock orthogonality.
We first observe that, for every entire function $g$,
\begin{equation}\label{eq:origin-test}
 g\in\mathcal N(I)\quad\Longleftrightarrow\quad
 (p^*(\partial)g)(0)=0\quad\text{for every }p\in I.
\end{equation}
The forward implication is immediate. For the converse, testing
$z^\beta p\in I$ gives
$\partial^\beta(p^*(\partial)g)(0)=0$ for every $\beta$, so the Taylor coefficients of $p^*(\partial)g$ must all be zero, and we get that $p^*(\partial)g$ vanishes. For $g\in\cF_d$ the monomial norms give, for $g(z)=\sum_\alpha g_\alpha z^\alpha$,
\[
 (p^*(\partial)g)(0)=\sum_\alpha\overline{p_\alpha}\,\alpha!\,g_\alpha
 =\langle g,p\rangle_{\cF_d}.
\]
Thus, $\mathcal N(I)\cap\cF_d=I^\perp$ in Fock space. The formal power series satisfying these equations form the inverse system of the ideal obtained by conjugating the coefficients of $I$, whose description goes back to Macaulay \cite[Ch.~IV, \S\S59--61]{macaulay1916}. Here, the series must define entire functions. Also, note that $\mathcal A(I)$ is an ideal of $\mathcal O(\C^d)$. Indeed, if $g$ is a limit of polynomials from $I$ and $h$ is entire, Taylor polynomials of $h$ multiplied by the approximants of $g$ give approximants of $hg$ in $I$.

The compact-open topology also controls differentiation. If
$g_\nu\to g$ uniformly on compact sets, then for each compact $E$,
$t>0$, and $\beta\in\Nz^d$, Cauchy's formula gives
\[
 \sup_{z\in E}|\partial^\beta(g_\nu-g)(z)|
 \leq \frac{\beta!}{t^{|\beta|}} \cdot
 \sup_{\substack{z\in E\\ |w_j-z_j|\leq t \\  \text{for} \  1\leq j\leq d}}
 |(g_\nu-g)(w)|\longrightarrow0 .
\]
By Montel's theorem, a sequence bounded uniformly on each compact subset of $\C^d$ has a subsequence converging in the compact-open topology. We will apply this compactness after moving the origin to a point where a function has large Gaussian-weighted absolute value. The translation must account for the change in the Gaussian as well as the change of ideal.

\begin{definition}\label{def:weighted-translation}
 For an entire function $g$, define its \emph{weighted absolute value} by $A_g(z)=e^{-|z|^2/2}|g(z)|$.
 For $b\in\C^d$ the \emph{weighted translate} of $g$ by $b$ is the entire
 function
 \[
 (D_b g)(z)=e^{-\langle z,b\rangle-|b|^2/2}g(z+b).
 \]
 The translate of an affine ideal $I$ by $b$ is
 $I_b=\{p(b+\cdot):p\in I\}$.
\end{definition}

The operator $D_b$ is what we call the \emph{Fock-space weighted translation} with translation parameter $-b$, or sometimes just \emph{translation} when the context is clear. Its exponential factor is chosen so that
\begin{equation}\label{eq:weighted-translation}
 A_{D_b g}(z)=e^{-|z|^2/2}|D_b g(z)|
             =e^{-|z+b|^2/2}|g(z+b)|=A_g(z+b).
\end{equation}
The modulus identity follows by expanding the Euclidean squared norm $|z+b|^2$ according to the inner product on $\mathbb{C}^d$. For $p\in I$, we have that $D_b p$ is a multiple by an entire function of $p(b+\cdot)\in I_b$ and must therefore belong to $\mathcal A(I_b)$.
The map $D_b$ is continuous for uniform convergence on compact sets.
Applying it to polynomial approximants gives
$D_b\mathcal A(I)\subseteq\mathcal A(I_b)$, and $D_{-b}D_b=1$ yields
equality. For the differential equations condition, the product rule yields
\begin{equation}\label{eq:translation-product-rules}
 D_b(z_j g)=(z_j+b_j)D_b g,\qquad
 D_b(\partial_j g)=(\partial_j+\overline{b_j})D_b g .
\end{equation}
Conjugating the coefficients of $p(b+z)$ replaces the shift $b$ by
$\overline b$, so $(p(b+\cdot))^*(\partial)=p^*(\partial+\overline b)$.
Iterating the second product-rule identity therefore gives
\[
 (p(b+\cdot))^*(\partial)D_b g=D_b(p^*(\partial)g).
\]
Applying this identity with $-b$ then gives $D_b\mathcal N(I)=\mathcal N(I_b)$.

These spaces behave naturally under a unitary change of coordinates $z\mapsto Wz$. Under composition with $W$, elements of $\mathcal A(I)$
are approximated by the transformed ideal polynomials. For the adjoint
equations the chain rule gives, for each coordinate index $a$,
\[
 \sum_j\overline{W_{aj}}\,\partial_j(g\circ W)
 =\sum_{j,k}\overline{W_{aj}}W_{kj}\,(\partial_k g)\circ W
 =(\partial_a g)\circ W .
\]
Iteration gives $(p\circ W)^*(\partial)(g\circ W)=(p^*(\partial)g)\circ W$.
The Gaussian weight is preserved since $|Wz|=|z|$.

\subsection{Homogenization and the Hilbert function}\label{subsec:filtered-homogenization}

To take a limit of the translated ideals $I_b$, we use the filtration
by total degree. Translation preserves each degree bound, and
homogenization records each filtered piece as a subspace of one
finite-dimensional homogeneous space. Translation keeps its dimension
fixed. We can therefore take limits of the subspaces degree by degree,
and the compatibility of multiplication will make their limit an ideal.

We introduce one additional variable $x_0$. In the remainder of this
section and in Section~\ref{sec:affine-limits}, $S=\C[x_0,z_1,\ldots,z_d]$
denotes the homogenizing ring, $S_m$ its homogeneous part of degree
$m$, and $R_{d,\leq m}$ the polynomials in $R_d$ of degree at most $m$.

\begin{definition}\label{def:homogenization}
 For $p=\sum_{|\alpha|\leq m}p_\alpha z^\alpha\in R_{d,\leq m}$,
 its \emph{homogenization to degree $m$} is
 \[
 p^{[m]}(x_0,z)=\sum_{|\alpha|\leq m}
                   p_\alpha x_0^{m-|\alpha|}z^\alpha .
 \]
 For an affine ideal $I\subseteq R_d$, its homogenization $I^h$ is the
 graded subspace of $S$ whose degree-$m$ piece is
 $(I^h)_m=\{p^{[m]}:p\in I\cap R_{d,\leq m}\}$.
 For a homogeneous ideal $H\subseteq S$, its \emph{dehomogenization}
 $H|_{x_0=1}$ is its image under $F\mapsto F(1,\cdot)$.
\end{definition}

Multiplication by $x_0$ or $z_j$ sends $p^{[m]}$ to $p^{[m+1]}$ or
$(z_jp)^{[m+1]}$, respectively. Thus, $I^h$ is a homogeneous ideal.
Evaluation at $x_0=1$ is a vector space isomorphism $S_m\to R_{d,\leq m}$
whose inverse is $p\mapsto p^{[m]}$. Moreover, every polynomial
in $H|_{x_0=1}$ has a lift in a single $H_m=H\cap S_m$. Indeed, if
$\sum_{j=0}^m F_j$ is any lift, with $F_j\in H_j$, then
$\sum_{j=0}^m x_0^{m-j}F_j$ is a lift in $H_m$.

Within each $S_m$ we compare subspaces through their orthogonal projections. The particular inner product used for this comparison is immaterial in finite dimension. Prescribing the dimensions of the quotient spaces keeps the ranks of these projections fixed.

\begin{definition}\label{def:degreewise-convergence}
 The \emph{Hilbert function} of a homogeneous quotient $S/H$ is
 $m\mapsto\dim_\C(S_m/H_m)$. Give $S_m$ the inner product for
 which its monomials are orthonormal. For a homogeneous ideal $H \subseteq S$ and a sequence $\{H_\nu\}_{\nu \in \mathbb{N}}$ of homogeneous ideals, all with the same quotient Hilbert function, we
 say that $H_\nu$ \emph{converges to $H$ in each degree} if
 \[
  \|P_{(H_\nu)_m}-P_{H_m}\|\longrightarrow0
  \qquad\text{for every }m\geq0 .
 \]
\end{definition}

By Definition~\ref{def:homogenization}, evaluation at $x_0=1$
identifies $S_m/(I^h)_m$ with $R_{d,\leq m}/(I\cap R_{d,\leq m})$.
Translation by $b$ carries $I^h$ onto $(I_b)^h$ via the graded
automorphism $F(x_0,z)\mapsto F(x_0,z+bx_0)$. Hence, the quotients by
these translated homogenizations all have the same Hilbert function.

For $0\leq h\leq\dim S_m$, the Grassmannian $\operatorname{Gr}(h,S_m)$
of $h$-dimensional subspaces of $S_m$ is identified with the rank-$h$
orthogonal projections. The conditions $P=P^*=P^2$ and
$\operatorname{Tr}P=h$ give a closed, bounded set in finite
dimension. These Grassmannians are therefore compact in the topology
of Definition~\ref{def:degreewise-convergence}. We use their compactness
to extract a limiting ideal.\footnote{Haiman and Sturmfels construct
a projective parameter space for ideals with a fixed Hilbert function
\cite[Theorem~1.1 and Corollary~1.2]{haimansturmfels2004}.}

\begin{lemma}\label{lem:filtered-compactness}
Every sequence of homogeneous ideals in $S$ with a common Hilbert function
has a subsequence converging in each degree to a homogeneous ideal
with the same Hilbert function.
\end{lemma}
\begin{proof}
 First extract a subsequence whose projections converge in degree zero,
 then a subsequence of it with convergence in degree one, and so on.
 The diagonal subsequence converges in every fixed degree to projections
 $P_{H_m}$ with the prescribed ranks. It remains to check the ideal
 property. For $F\in H_m$, let
 $F_\nu=P_{(H_\nu)_m}F$. Then, $F_\nu\to F$. For every variable $x$
 of $S$, $xF_\nu\to xF$ and $P_{(H_\nu)_{m+1}}xF_\nu=xF_\nu$.
 Taking limits gives $P_{H_{m+1}}xF=xF$, so $xH_m\subseteq H_{m+1}$.
 Therefore, $\bigoplus_{m\geq0}H_m$ is a homogeneous ideal.
\end{proof}

\begin{example}\label{ex:homogeneous-limit}
A homogeneous limit can differ from the homogenization of its
dehomogenization. The ideal $I_\nu=(z+\nu)$ in one variable has homogenization
$(x_0+z/\nu)$, whose degreewise limit is $H=(x_0)$. Its affine image is $L=(1)$, so $L^h=S\ne H$.
\end{example}

\subsection{Coefficient norms and leading monomials}\label{subsec:coefficient-norms}

Convergence in each fixed degree gives no estimate for a polynomial
whose degree increases along the sequence. We choose norms for which
an estimate on finitely many generators also controls all their monomial multiples, with no dependence on the degree. These norms will
also let us pass from convergence of polynomial coefficients to locally uniform convergence of entire functions.

\begin{definition}\label{def:coefficient-norm}
For positive radii $r=(r_1,\ldots,r_d)$, let $\Delta_r$ be the closed
polydisc $|z_j|\le r_j$, and write $\Delta_R$ when all the radii are equal to $R$. For an entire function $f(z)=\sum_\alpha f_\alpha z^\alpha$, put $r^\alpha=\prod_j r_j^{\alpha_j}$. The \emph{coefficient norm} of $f$ at $r$ is $\|f\|_r=\sum_\alpha |f_\alpha|r^\alpha$.
We write $\Pi_N f=\sum_{|\alpha|\le N}f_\alpha z^\alpha$ for its
Taylor polynomial truncation up to total degree $N$.
\end{definition}

Cauchy's inequalities give
$|f_\alpha|r^\alpha\le 2^{-|\alpha|}\sup_{z\in \Delta_{2r}}|f(z)|$.
Summing and using $\sum_\alpha 2^{-|\alpha|}=2^d$ proves both finiteness
of the coefficient norm and the comparison between the norms
\begin{equation}\label{eq:coefficient-norm-comparison}
 \sup_{\Delta_r}|f|\ \le\ \|f\|_r\ \le\ 2^d\sup_{\Delta_{2r}}|f|.
\end{equation}
Thus, the norms $\|f\|_r$ induce the topology of uniform convergence on compact sets. For homogeneous elements in $S_m$ we use the same notation, giving $x_0$ radius one. In a monomial $x_0^{m-|\alpha|}z^\alpha$, the exponent of $x_0$ is determined by $\alpha$ and $m$. Hence, setting $x_0=1$ combines no two distinct terms, and we get that for $F\in S_m$,
\begin{equation}\label{eq:homogeneous-coefficient-norm}
 \|F\|_r=\|F(1,\cdot)\|_r.
\end{equation}
Also, note that $\|uF\|_r=\|u\|_r\|F\|_r$ for every monomial $u$.

We next choose leading monomials so that they dominate the other terms
of finitely many forms when the radii are sufficiently large. Fix positive numbers $a_1,\ldots,a_d$ linearly independent over $\mathbb Q$. Give $x_0^{m-|\alpha|}z^\alpha\in S_m$ the weight which is the dot product $a \cdot \alpha$. Thus, any two weights in a fixed total degree are distinct. The \emph{leading monomial} of a nonzero homogeneous polynomial is the monomial of greatest weight with nonzero coefficient. Ordering all monomials first by total degree and then by weight gives a \emph{monomial} order on monomials $S$. Multiplication preserves this order, and each bounded total degree contains only finitely many monomials.

For a homogeneous ideal $H\ne0$, its possible leading monomials are
closed under multiplication by monomials. Their span is the
\emph{initial ideal} $\operatorname{in}(H)$. By the Hilbert basis theorem, it has
finitely many minimal monomial generators $M_1,\ldots,M_s$.
Choose homogeneous $G_i\in H$ whose leading term is $M_i$ with
coefficient one. Such a family is then a homogeneous Gr\"obner basis
\cite{eisenbud}.

Let $E_m=(\operatorname{in}(H))_m$ be the span of the degree-$m$ monomials in $\operatorname{in}(H)$. By \cite[Theorem~15.26, p.~356]{eisenbud}, the quotients $S/H$ and $S/\operatorname{in}(H)$ have the same Hilbert function. Consequently,
\[
 \dim E_m
 =\dim S_m-\dim(S/\operatorname{in}(H))_m
 =\dim S_m-\dim(S/H)_m
 =\dim H_m.
\] The finite family $G_1,\ldots,G_s$ and this dimension identity are the facts about Gr\"obner bases that we will use, and no division algorithm for Gr\"obner bases is needed below.

\section{Limits of affine ideals}\label{sec:affine-limits}

This section supplies the two properties of limiting ideals needed
in the proof of our main result. First, both approximation by ideal polynomials and the differential equations for orthogonality survive passage to a limit. Second, a limit of translates of a fixed ideal, with centers tending to infinity, is preserved by a nonzero translation. This extra symmetry will allow us to remove one variable by induction in Section~\ref{sec:affine-uniqueness}.

Throughout the Subsection \ref{subsec:passing-to-limit}, let $H_\nu,H\subseteq S$ be homogeneous ideals with the same quotient Hilbert function, and suppose that $H_\nu$ converges to $H$ in each degree. Put $I_\nu=H_\nu|_{x_0=1}$ and $L=H|_{x_0=1}$. We use the actual homogeneous ideals $H_\nu$ and $H$, so that $H$ need not coincide with the homogenization of its dehomogenization.

\subsection{Passing to the limit}\label{subsec:passing-to-limit}

The two conditions on entire functions present different problems.
To pass the differential equations to a limit, we can test one
polynomial of bounded degree at a time. For approximation by ideal
polynomials, the degrees of the approximants may tend to infinity.
We therefore need an estimate that works in every degree at once.
The idea is to compare finitely many polynomials in $H_\nu$ with fixed
monomials and then use that same comparison for all of their monomial
multiples.

For a fixed ideal, the use of weighted coefficient norms to make the leading monomials dominate goes back to Djakov and Mitiagin \cite[equations (1.12)--(1.15)]{djakovmitiagin}. See also \cite[\S3]{astw1996} and \cite[Proposition~3.2 and Lemma~3.4]{shalitdivision}. Here the polynomials used for division vary with $\nu$. We choose the radii using the limiting ideal and show that the resulting estimate holds in every degree for all sufficiently large $\nu$.

Suppose $H\ne0$, and choose $M_i,G_i,E_m$ as in Subsection~\ref{subsec:coefficient-norms}. Write $d_i=\deg G_i$ and put
$G_{i,\nu}=P_{(H_\nu)_{d_i}}G_i$. Convergence in degree $d_i$ gives
$G_{i,\nu}\to G_i$ coefficientwise. Assign each monomial
$M\in\operatorname{in}(H)$ to the first $M_i$ that divides it.
This assignment is fixed in all degrees and for all $\nu$.
We will say that a polynomial coefficient $C_i$, for $i=1,...,s$, \emph{respects the assignment} if every monomial $u$ occurring in $C_i$ has $uM_i$ assigned to $M_i$. If $u$ occurs in $C_i$ and $v$ occurs in $C_j$, then $uM_i = vM_j$ forces $i= j$, since this monomial has a unique
assigned divisor. Cancellation then gives $u=v$. Consequently, no monomial occurs in two different products among $M_1C_1,\ldots,M_sC_s$. The assignment also defines a linear map
\[
 A_{\nu,m}:E_m\longrightarrow(H_\nu)_m,\qquad
 A_{\nu,m}(uM_i)=uG_{i,\nu},
\]
where $M_i$ is the assigned divisor of the basis monomial $uM_i$.
Its range lies in $(H_\nu)_m$ because $H_\nu$ is an ideal. We will
choose coefficient norms in which this map is close to the inclusion
of $E_m$ into $S_m$, with an error independent of $m$.

\begin{lemma}\label{lem:assigned-monomials}
For each $R\ge1$ there are radii $r_j\ge R$ and an index $\nu_0$
such that the following holds for every $\nu\ge\nu_0$ and every
$m\ge0$. Each $F\in(H_\nu)_m$ has a unique expansion
\[
 F=\sum_{i=1}^s G_{i,\nu}C_i
\]
with homogeneous $C_i$ of degree $m-d_i$ respecting the assignment,
where $C_i=0$ if $m<d_i$. Moreover,
\begin{equation}\label{eq:assigned-monomials}
 \sum_{i=1}^s\|M_i\|_r\|C_i\|_r\le 2\|F\|_r .
\end{equation}
\end{lemma}
\begin{proof}
Put $r_j=t^{a_j}$ with $t>1$, using the weights chosen in
Subsection~\ref{subsec:coefficient-norms}. Every monomial of $G_i-M_i$
has strictly smaller weight than $M_i$, so
$\|G_i-M_i\|_r/\|M_i\|_r\to0$ as $t\to\infty$.
Choose $t$ large enough that all $r_j\ge R$ and the ratios $\|G_i-M_i\|_r/\|M_i\|_r$ are at most $1/4$. Fix these radii $r_j$. Convergence of coefficients of the finitely many $G_{i,\nu}$ then gives $\nu_0$ such that for all $i=1,\ldots,s$ and $\nu\ge\nu_0$,
\[
 \|G_{i,\nu}-M_i\|_r
 \le\|G_{i,\nu}-G_i\|_r+\|G_i-M_i\|_r
 \le\tfrac12\|M_i\|_r.
\]
In particular, $t$ and $\nu_0$ have been chosen without reference
to the degree $m$.

Fix $\nu\ge\nu_0$ and $m\ge0$. Exact scaling of the coefficient
norm under monomial multiplication gives
\[
 \|A_{\nu,m}(uM_i)-uM_i\|_r
 =\|u\|_r\|G_{i,\nu}-M_i\|_r\leq \tfrac12 \|u\|_r \|M_i\|_r =
 \tfrac12\|uM_i\|_r.
\]
For $v=\sum_M v_M M\in E_m$, where $M$ are monomials in $E_m$ with an assigned monomial among $M_1,...,M_s$, we therefore have
\[
 \|A_{\nu,m}v-v\|_r
 \le\sum_M |v_M|\,\|A_{\nu,m}M-M\|_r \leq \tfrac12 \sum_M |v_M| \|M\|_r
  = \tfrac12\|v\|_r.
\]
Hence, by reverse triangle inequality, $\|A_{\nu,m}v\|_r\ge\tfrac12\|v\|_r$ and we get that $A_{\nu,m}$ is injective. The common Hilbert function now gives
\[
 \dim(H_\nu)_m=\dim H_m=\dim E_m,
\]
so we deduce that $A_{\nu,m}$ is onto as well. This is where we use equality of dimensions in every degree, including degrees that may increase with $\nu$.

Given $F\in(H_\nu)_m$, group the monomials of $v=A_{\nu,m}^{-1}F$ according to their assigned divisors. This gives $v=\sum_iM_iC_i$, with each $C_i$ respecting the assignment, and an application of $A_{\nu,m}$ to this gives $F=\sum_iG_{i,\nu}C_i$. Conversely, any such expansion determines a preimage $v=\sum_iM_iC_i$ of $F$, which is unique by injectivity
of $A_{\nu,m}$. Writing $v=\sum_M a_M M$, the assignment then determines each polynomial coefficient uniquely by setting
\[
 C_i=\sum_{\substack{M\text{ assigned}\\\text{to }M_i}}
          a_M\,\frac{M}{M_i}.
\]
Since each monomial of $v$ occurs in exactly one of the products $M_iC_i$, by the definition of the coefficient norm $\| \cdot \|_r$ gives, together with the inequality $\|A_{\nu,m}v\|_r\ge\tfrac12\|v\|_r$, that
\[
 \sum_i\|M_i\|_r\|C_i\|_r = \sum_i\|M_iC_i\|_r =  \|v\|_r \le 2\|F\|_r. \qedhere
\] 
\end{proof}

The map $A_{\nu,m}$ does not depend on the radii. Consequently,
whenever the lemma applies with two choices of radii, it gives the
same coefficients. This will let us estimate one sequence of coefficients on successively larger polydiscs.

\begin{example}\label{ex:hilbert-function}
The common Hilbert function cannot be replaced by convergence of the pieces in each fixed degree alone. Let $H_\nu=((x_0+z/\nu^2)^\nu)$ in $\C[x_0,z]$. Its piece in each
fixed degree is eventually zero. Thus, each fixed-degree piece converges to the corresponding piece of the zero ideal. Nevertheless, the dehomogenized generators $(1+z/\nu^2)^\nu$ converge to $1$ uniformly on compact sets.
Thus, the limit of these ideal polynomials does not belong to the compact-open closure of the limiting ideal.
\end{example}

We can now pass both conditions to the limit. The proof of the approximation assertion will also produce a representation of a limit with entire coefficients, which we will use in the next section.

\begin{lemma}\label{lem:affine-limits}
Let $H_\nu,H,I_\nu,L$ be as above, and $G_i$ the chosen Gr\"obner basis for $H$, and suppose that entire functions
$g_\nu$ converge to $g$ uniformly on compact subsets of $\C^d$. Then, \[
 \begin{aligned}
 g_\nu\in\mathcal N(I_\nu)\text{ for all }\nu
 &\ \Longrightarrow\ g\in\mathcal N(L),\\
 g_\nu\in\mathcal A(I_\nu)\text{ for all }\nu
 &\ \Longrightarrow\ g\in\mathcal A(L).
 \end{aligned}
\]
In the second case, if $H\ne0$, there are entire functions $c_i$ such
that $g=\sum_iG_i(1,\cdot)c_i$.
\end{lemma}
\begin{proof}
Suppose first that $g_\nu\in\mathcal N(I_\nu)$. Fix $p\in L$ and
choose a homogeneous lift $F\in H_m$ with $F(1,\cdot)=p$, as in
Subsection~\ref{subsec:filtered-homogenization}. The polynomials
\[
 p_\nu=(P_{(H_\nu)_m}F)(1,\cdot)\in I_\nu
\]
have degree at most $m$ and converge coefficientwise to $p$.
Cauchy's formula gives local uniform convergence of all derivatives
of $g_\nu$ up to order $m$. Consequently,
\[
 p_\nu^*(\partial)g_\nu\longrightarrow p^*(\partial)g
 \quad\text{locally uniformly}.
\]
The functions on the left are identically zero, so
$p^*(\partial)g=0$. Since $p\in L$ was arbitrary, this proves
$g\in\mathcal N(L)$.

Now suppose $g_\nu\in\mathcal A(I_\nu)$. Choose a polynomial $q_\nu\in I_\nu$ within $1/\nu$ of $g_\nu$ uniformly on $\Delta_\nu$. Every fixed
compact set is eventually contained in $\Delta_\nu$ for sufficiently large $\nu$, so
$q_\nu\to g$ locally uniformly. If $H=0$, equality of the Hilbert
functions gives $(H_\nu)_m=0$ for all $\nu,m$, and hence $g=0$.
We may therefore assume $H\ne0$.

For each $\nu$, fix a homogeneous lift
$F_\nu\in(H_\nu)_{m_\nu}$ of $q_\nu$. The degrees $m_\nu$ may be
unbounded in $\nu$. Apply Lemma~\ref{lem:assigned-monomials} with $R=1$.
After discarding finitely many indices $\nu$, we obtain the unique
expansions $F_\nu=\sum_iG_{i,\nu}C_{i,\nu}$. Put
$c_{i,\nu}=C_{i,\nu}(1,\cdot)$, so that
\[
 q_\nu=\sum_iG_{i,\nu}(1,\cdot)c_{i,\nu}.
\]
We claim that each sequence $(c_{i,\nu})_\nu$ is bounded on every
compact set.

To prove the claim, fix $R\ge1$ and apply Lemma~\ref{lem:assigned-monomials} again with this $R$ to obtain radii $r_j\ge R$ and a tail of indices of $\nu$ on which the estimate of the lemma holds. The lifts and the assignment are unchanged. The coefficients for $\nu$ on this
tail are therefore the ones already chosen. Since $\|M_i\|_r\ge1$, equations \eqref{eq:assigned-monomials}, \eqref{eq:homogeneous-coefficient-norm}, and \eqref{eq:coefficient-norm-comparison} yield
\[
 \sup_{\Delta_R}|c_{i,\nu}| \le \| c_{i,\nu} \|_r = \|C_{i,\nu}\|_r\le2\|F_\nu\|_r
 =2\|q_\nu\|_r\le2^{d+1}\sup_{\Delta_{2r}}|q_\nu|.
\]
Here, $r$ is fixed before $\nu$ varies, so local uniform convergence
of $q_\nu$ bounds the last expression. The finitely many excluded
coefficients are polynomials and are bounded on $\Delta_R$ as well.
This proves the claim for the same coefficient sequences on every
polydisc.

By Montel's theorem, a subsequence of the finite tuple
$(c_{1,\nu},\ldots,c_{s,\nu})$ converges locally uniformly to entire
functions $(c_1,\ldots,c_s)$. Passing to the limit in the displayed
representation of $q_\nu$ gives $g=\sum_iG_i(1,\cdot)c_i$.
Finally, the polynomials $\sum_iG_i(1,\cdot)\Pi_Nc_i$ belong to
$L$ and converge locally uniformly to $g$, and it follows that $g\in\mathcal A(L)$.
\end{proof}

For a nonzero fixed ideal $I$, apply the lemma to the constant
sequences $H_\nu=H=I^h$ and $g_\nu=g$. It writes every
$g\in\mathcal A(I)$ as a finite sum $g=\sum_i p_i c_i$ with
$p_i=G_i(1,\cdot)\in I$ and $c_i$ entire. This is the representation
needed in Section~\ref{sec:affine-uniqueness}.

\subsection{A translation direction at infinity}

We next consider translates of a fixed ideal with centers tending
to infinity. The aim is to find a nonzero translation that maps the
limiting affine ideal into itself. We will compare two translated
ideals whose centers are one unit apart and whose homogeneous pieces
are both close to the limit. Their displacements have a convergent
subsequence, and the limiting displacement gives the desired
translation.

We will compare projections using the Hilbert--Schmidt norm, which is the square root of the sum of the squared absolute values of the matrix entries in the orthonormal monomial basis. We also use one fact about real polynomial inequalities. A \emph{semialgebraic set} is obtained by finitely many unions, intersections, and complements
from sets defined by real polynomial inequalities and equalities.
We will make use of the fact that every such set has finitely many connected components \cite[Theorem~2.4.5]{bochnakcosteroy1998}, and will identify
$\C^d$ with $\mathbb R^{2d}$ when applying this fact.

\begin{lemma}\label{lem:translation-line}
Let $I\subseteq R_d$ be an affine ideal and suppose that
$|b_\nu|\to\infty$ while $(I_{b_\nu})^h\to H$ in each degree.
There is a unit vector $e\in\C^d$ such that translation by $e$ maps
$L=H|_{x_0=1}$ into itself.
\end{lemma}

\begin{proof}
For $b\in\C^d$, write $\tau_bF(x_0,z)=F(x_0,z+bx_0)$ for the graded translation on $S$. Then, we have that $\tau_b\tau_c=\tau_{b+c}$ and $\tau_b^{-1}=\tau_{-b}$. If $H=0$, then $L=0$, and any unit vector works. Otherwise, choose nonzero homogeneous generators $F_1,\ldots,F_s$ of $H$
and put $D=\max_i\deg F_i$. Since
$x_0^{D-\deg F_i}F_i\in H_D$ has dehomogenization $F_i(1,\cdot)$,
the image of $H_D$ under evaluation at $x_0=1$ generates $L$.
It therefore suffices to find a unit vector $e$ such that
$\tau_eH_D\subseteq H_D$.

Put $E_0=(I^h)_D$ and $E(b)=\tau_bE_0=((I_b)^h)_D$.
Let $P(b)$ and $P$ be the orthogonal projections of $S_D$ onto
$E(b)$ and $H_D$, respectively, for the coefficient inner product.
By hypothesis, $P(b_\nu)\to P$ in the operator norm. We will find
pairs of centers one unit apart whose projections both tend to
$P$, and then take a limit of their displacements.

To describe these projections explicitly, fix an ordering of the
orthonormal monomial basis
$m_\alpha=x_0^{D-|\alpha|}z^\alpha$, indexed by
$\mathcal M_D=\{\alpha\in\mathbb Z_{\ge0}^d:|\alpha|\le D\}$.
Write $N=|\mathcal M_D|=\dim S_D$, and choose a fixed basis
$Q_1,\ldots,Q_h$ of $E_0$. Write
$Q_k=\sum_{\alpha\in\mathcal M_D}q_{\alpha k}m_\alpha$.
Define $V(b)$ to be the $N\times h$ matrix whose $k$th column
is the coefficient vector of $\tau_bQ_k$ in this monomial basis.
The binomial theorem gives
\[
 V_{\beta k}(b)
 =
 \sum_{\substack{\alpha\in\mathcal M_D\\\alpha\ge\beta}}
 q_{\alpha k}\binom{\alpha}{\beta}b^{\alpha-\beta},
 \qquad
 \binom{\alpha}{\beta}
 =\prod_{j=1}^d\binom{\alpha_j}{\beta_j},
\]
where $\alpha\ge\beta$ means coordinatewise inequality.
Thus, the entries of $V(b)$ are polynomial in $b$, and its
columns correspond to a basis of $E(b)$ because $\tau_b$ is invertible.

The positive semi-definite Gram matrix $B(b)=V(b)^*V(b)$ is therefore invertible, and we have that $P(b)=V(b)B(b)^{-1}V(b)^*$. Put $\eta(b)=\det B(b)>0$ and $A(b)=V(b)\operatorname{adj}(B(b))V(b)^*$, where $\operatorname{adj}$ denotes the classical adjugate. Then, we have that $P(b)=A(b)/\eta(b)$. In the real coordinates
$b_j=u_j+iv_j$, the entries of $A(b)$ are polynomials and
$\eta(b)$ is a real polynomial.

For each positive integer $\ell$, let $U_\ell
 =\{b\in\C^d:\|P(b)-P\|_{\mathrm{HS}}<1/\ell\}$. Since $\eta(b)>0$, the defining inequality is equivalent to
\[
 \ell^2 \sum_{\alpha,\beta\in\mathcal M_D}
 \left|A_{\alpha\beta}(b)-\eta(b)P_{\alpha\beta}\right|^2
 <\eta(b)^2.
\]
This is a real polynomial inequality, since each squared absolute
value is the sum of the squares of its real and imaginary
parts. Hence, $U_\ell$ is a semialgebraic set. Moreover,
$\|T\|_{\mathrm{HS}}\le\sqrt N\,\|T\|$ on $S_D$, so
$P(b_\nu)\to P$ also in Hilbert--Schmidt norm. Thus, $U_\ell$
contains a tail of $(b_\nu)$ and is unbounded.

Since a semialgebraic set has finitely many connected components,
$U_\ell$ has an unbounded connected component $C_\ell$.
Fix $a_\ell\in C_\ell$. The continuous function
$b\mapsto|b-a_\ell|$ has connected, unbounded image on $C_\ell$,
and its image contains zero. This image is therefore an interval
containing $1$, so we can choose $c_\ell\in C_\ell$ with
$|c_\ell-a_\ell|=1$. Both points belong to $U_\ell$, and hence
$P(a_\ell)\to P$ and $P(c_\ell)\to P$.

Put $h_\ell=c_\ell-a_\ell$. After passing to a subsequence,
compactness of the unit sphere gives $h_\ell\to e$ with $|e|=1$.
Fix $F\in H_D$ and put $F_\ell=P(a_\ell)F$. Then,
$F_\ell\in E(a_\ell)$ and $F_\ell\to PF=F$. Since
$\tau_{h_\ell}\tau_{a_\ell}=\tau_{c_\ell}$, we have
\[
 \tau_{h_\ell}F_\ell\in
 \tau_{h_\ell}E(a_\ell)=E(c_\ell).
\]
The matrix entries of $\tau_b$ on $S_D$ are polynomial in $b$,
as explained above using the binomial expansion. Consequently,
$\tau_{h_\ell}\to\tau_e$ in operator norm, and therefore
$\tau_{h_\ell}F_\ell\to\tau_eF$. As $P(c_\ell)$ fixes
$\tau_{h_\ell}F_\ell$, passing to the limit gives
\[
 P\tau_eF
 =\lim_{\ell\to\infty}P(c_\ell)\tau_{h_\ell}F_\ell
 =\lim_{\ell\to\infty}\tau_{h_\ell}F_\ell
 =\tau_eF.
\]
Thus, $\tau_eF\in H_D$, proving $\tau_eH_D\subseteq H_D$.

Finally, write $p_i(z)=F_i(1,z)$. Applying this inclusion to
$x_0^{D-\deg F_i}F_i$ and evaluating at $x_0=1$ gives
$p_i(z+e)\in L$. Every $p\in L$ has an expression
$p=\sum_i a_ip_i$ with polynomial coefficients $a_i$, so
$p(z+e)=\sum_i a_i(z+e)p_i(z+e)\in L$.
Hence, translation by the unit vector $e$ maps $L$ into itself.
\end{proof}

The complex hyperplane perpendicular to $e$ contains the origin. Restriction
to it therefore retains the value of an entire function at zero and
any bound on its Gaussian-weighted absolute value. We must also
retain the two ideal conditions. The algebraic point is that a
polynomial ideal preserved by one nonzero translation is closed
under restriction to the perpendicular hyperplane. Polynomial
interpolation will express this restriction as a polynomial
combination of finitely many integer translates.

For $q\ge0$, the Lagrange polynomials for the points $0,\ldots,q$ are, for $0\le j\le q$,
\[
 \ell_{j,q}(s)=\prod_{\substack{0\le k\le q\\ k\ne j}}
                  \frac{s-k}{j-k},
\]
with the empty product equal to $1$. Lagrange interpolation gives
$F=\sum_{j=0}^q F(j)\ell_{j,q}$ for every polynomial $F$ of degree
at most $q$, also when its coefficients are polynomials in other variables.

\begin{lemma}\label{lem:transverse-restriction}
Suppose that translation by a unit vector $e$ maps an ideal
$L\subseteq R_d$ into itself. In unitary coordinates $(w,z')$ in
which $e=(1,0,\ldots,0)$, put $L'=L\cap\C[z']$. If
$G\in\mathcal A(L)\cap\mathcal N(L)$, then
$G(0,\cdot)\in\mathcal A(L')\cap\mathcal N(L')$.
\end{lemma}
\begin{proof}
Let $p\in L$ be nonzero and put $q=\deg p$. Apply interpolation
in the variable $s$ to $F(s)=p(w+s,z')$, whose coefficients belong
to $\C[w,z']$. Substituting $s=-w$ in the Lagrange interpolation polynomial identity gives
\[
 p(0,z')=\sum_{j=0}^q\ell_{j,q}(-w)\,p(w+j,z').
\]
Each translate $p(w+j,z')$ belongs to $L$ by repeated application of
the hypothesis, and each coefficient $\ell_{j,q}(-w)$ is a polynomial
in $w$. Thus, $p(0,z')\in L\cap\C[z']=L'$.
The same conclusion holds for $p=0$.

Restricting polynomial approximants from $L$ now proves
$G(0,\cdot)\in\mathcal A(L')$. For $p\in L'$, the operator
$p^*(\partial)$ involves only derivatives in $z'$, so it commutes
with restriction to the hyperplane $w=0$. Thus,
\[
 p^*(\partial_{z'})G(0,z')=(p^*(\partial)G)(0,z')=0,
\]
where the last equality uses $p\in L$. This proves
$G(0,\cdot)\in\mathcal N(L')$.
\end{proof}

We collect the limit statements in the form used in the next two
sections. The subsequence is chosen using only the centers. It therefore works for both conditions on entire functions, and remains suitable if we subsequently pass to a further subsequence to obtain convergence of functions.

\begin{lemma}\label{lem:translation-limits}
Let $I\subseteq R_d$ be an affine ideal and let
$(b_\nu)_{\nu \in \mathbb{N}}$ be a sequence in $\C^d$. There are a subsequence and an affine ideal $L$ such that, along this subsequence,
every locally uniform limit of functions
$f_\nu\in\mathcal A(I_{b_\nu})$ belongs to $\mathcal A(L)$,
and every locally uniform limit of functions
$g_\nu\in\mathcal N(I_{b_\nu})$ belongs to $\mathcal N(L)$.
If $|b_\nu|\to\infty$, translation by some unit vector
maps $L$ into itself.
\end{lemma}
\begin{proof}
Put $H_\nu=(I_{b_\nu})^h$. The graded automorphism
$F(x_0,z)\mapsto F(x_0,z+b_\nu x_0)$ carries $I^h$ onto
$H_\nu$, so the quotients $S/H_\nu$ have a common Hilbert
function. By Lemma~\ref{lem:filtered-compactness}, after
passing to a subsequence, the ideals $H_\nu$ converge
in each degree to a homogeneous ideal $H$.
Set $L=H|_{x_0=1}$.

Lemma~\ref{lem:affine-limits} gives both assertions about
locally uniform limits. The subsequence and $L$ were
chosen using only the translated ideals, so the same
choice works for both assertions and for every convergent
sequence of functions as in the statement.
If $|b_\nu|\to\infty$, Lemma~\ref{lem:translation-line}
gives a unit vector whose translation maps $L$ into itself.
\end{proof}

\section{Uniqueness for affine ideals}\label{sec:affine-uniqueness}

We now show that the two conditions preserved in the preceding
section are incompatible for a nonzero entire function with suitable
Gaussian growth. More precisely, if
$g\in\mathcal A(I)\cap\mathcal N(I)$ and
$A_g(z)=e^{-|z|^2/2}|g(z)|$ has polynomial growth, then $g=0$.
The application in Section~\ref{sec:singular-vectors} will produce a
function with bounded $A_g$. 

We say that $A_g$ has \emph{polynomial growth} if
$A_g(z)\le C(1+|z|)^m$ for some $C,m\ge0$, and that $A_g$ is
\emph{rapidly decreasing} if $(1+|z|)^kA_g(z)$ is bounded for every
integer $k\ge0$. The proof splits according to whether or not $A_g$ is
rapidly decreasing. In the case when $A_g$ is rapidly decreasing, a division theorem supplies a
representation of $g$ with entire coefficients that can be
paired with $g$. Polynomial truncations then reduce the pairing to Fock orthogonality. Otherwise, when $A_g$ is not rapidly decreasing, we obtain a nonzero limit of normalized translates with centers tending to infinity. A translation preserving
the limiting ideal allows us to restrict to a hyperplane and use
induction on the number of variables. This normalization may increase
the growth exponent, which is why the induction treats arbitrarily large
polynomial growth exponent.

For principal ideals, Fock approximation has a classical precedent.
Newman and Shapiro proved that every entire multiple of a nonzero
polynomial $p$ that belongs to $\cF_d$ lies in the Fock closure of
$pR_d$ \cite{newmanshapiro1968}. If $p$ is homogeneous, Shapiro's
extension of Fischer's decomposition gives uniqueness for entire
multiples satisfying $p^*(\partial)g=0$ without a growth assumption
\cite[Theorem~1]{shapiro1989}. Meril and Yger proved uniqueness
among entire functions for a single nonzero polynomial of degree at most
two, and for systems homogeneous with respect to common positive integer
weights whose generators form a regular sequence in every order
\cite[Theorem~3 and Proposition~8]{merilyger1992}.

Write $d\gamma(z)=e^{-|z|^2}d\mu(z)$ for the Fock measure, and
$|h|^2=\sum_i|h_i|^2$ for a vector of functions. For a polynomial
row $P=(p_1,\ldots,p_s)$, we similarly put $|P|^2=\sum_i|p_i|^2$.
The representation
supplied by Lemma~\ref{lem:affine-limits} has entire coefficients,
but gives no control of their growth at infinity. We need a
representation whose coefficients have a weighted Gaussian norm.
The division estimate below in Lemma \ref{lem:gaussian-division} permits a polynomial loss in the Gaussian norm of the coefficients. Rapid decrease of $A_g$ will absorb this loss when we form the pairing. The estimate uses
\emph{plurisubharmonic} weights. For a smooth real function $\phi$,
this means that its complex Hessian
$(\partial^2\phi/\partial z_j\partial\overline z_k)_{j,k}$ is
positive semidefinite. Here, writing $z_j=x_j+iy_j$, we use
$\partial/\partial z_j=(\partial/\partial x_j-i\partial/\partial y_j)/2$
and $\partial/\partial\overline z_j=(\partial/\partial x_j+i\partial/\partial y_j)/2$.
For holomorphic functions, the first derivative agrees with $\partial_j$
used earlier. All weights below are smooth. 

In what follows, we will make use of H\"ormander's division theorem
\cite[equation (7.6.16), Theorem~7.6.11 and equation (7.6.19),
pp.~219--222]{hormander}, which we now describe. For a fixed row of polynomials $P$ and a vector of entire functions $u$, we denote $Pu=\sum_i p_i u_i$.  The theorem requires two plurisubharmonic functions $\phi$ and $-\log\rho$ with $0<\rho\le1$, satisfying
\[
 \begin{aligned}
 \rho(z+\xi)&\le2\rho(z)
 &&\text{if }|\mathrm{Re}\xi_j|,|\mathrm{Im}\xi_j|\le1
   \text{ for every }j,\\
 |\phi(z+\xi)-\phi(z)|&\le C_0
 &&\text{if }|\xi_j|<\rho(z)
   \text{ for every }j.
 \end{aligned}
\]
The theorem produces an integer $N\ge0$
such that, whenever $u$ is a vector of entire functions with
$\int_{\C^d}|Pu|^2e^{-\phi}\,d\mu<\infty$, there is a vector of entire functions $h$ satisfying $Ph=Pu$ and
\[
 \int_{\C^d}|h|^2e^{-\phi}\rho^N(1+|z|^2)^{-N}\,d\mu
 \le C\int_{\C^d}|Pu|^2e^{-\phi}\,d\mu.
\]
$N$ only depends on $P$, and $C$ depends only on $P$
and $C_0$. No integrability assumption is imposed on the
original entries of $u$.

\begin{lemma}\label{lem:gaussian-division}
For each fixed polynomial row $P=(p_1,\ldots,p_s)$, there are
constants $\kappa\ge0$ and $C_P\ge0$ such that the following holds. If $u$ is
a vector of entire functions such that $Pu \in\cF_d$, then there is a vector of entire functions $h$ with $Ph=Pu$ and
\begin{equation}\label{eq:hormander-division}
 \int_{\C^d}|h(z)|^2(1+|z|^2)^{-\kappa}\,d\gamma(z)
 \le C_P\int_{\C^d}|Pu(z)|^2\,d\gamma(z).
\end{equation}
\end{lemma}
\begin{proof}
We apply H\"ormander's division theorem with
$\phi(z)=|z|^2$ and
$\rho(z)=(a_0+|z|^2)^{-1/2}$, where $a_0=1+2d$.
Denote $A(z)=a_0+|z|^2$. It is easy to prove that $\phi$ and $-\log \rho$ are plurisubharmonic, and we also have that
$0<\rho\le a_0^{-1/2}\le1$.

To verify the first translation condition, suppose that
$|\mathrm{Re} \xi_j|,|\mathrm{Im}\xi_j|\le1$ for every $j$. Then
$|\xi|\le\sqrt{2d}\le\sqrt{a_0}$, and the triangle inequality
applied to the vector $(\sqrt{a_0},z)$ gives
\[
 \sqrt{a_0+|z|^2}
 \le \sqrt{a_0+|z+\xi|^2}+|\xi|
 \le 2\sqrt{a_0+|z+\xi|^2},
\]
which establishes $\rho(z+\xi)\le2\rho(z)$.
For the second condition, if $|\xi_j|<\rho(z)$ for every $j$,
then $|\xi|\le\sqrt d\,A(z)^{-1/2}$, and hence
\[
 \begin{aligned}
 |\phi(z+\xi)-\phi(z)|
 \le 2|z||\xi|+|\xi|^2 \le \frac{2\sqrt d\,|z|}{\sqrt{A(z)}}+\frac d{A(z)}
 \le 2\sqrt d+\frac d{a_0}.
 \end{aligned}
\]
The theorem therefore applies with
$C_0=2\sqrt d+d/a_0$.

It remains to compare the chosen weight with the one in the statement.
Since $a_0\ge1$, we have
$A(z)\le a_0(1+|z|^2)$, so
\[
 \rho(z)^N(1+|z|^2)^{-N}
 =A(z)^{-N/2}(1+|z|^2)^{-N}
 \ge a_0^{-N/2}(1+|z|^2)^{-3N/2}.
\]
Using $d\gamma=e^{-|z|^2}\,d\mu$, we consequently obtain
\[
 \begin{aligned}
 \int_{\C^d}|h|^2(1+|z|^2)^{-3N/2}\,d\gamma \le a_0^{N/2}
       \int_{\C^d}|h|^2\rho^N(1+|z|^2)^{-N}\,d\gamma \le C a_0^{N/2}\int_{\C^d}|Pu|^2\,d\gamma.
 \end{aligned}
\]
This proves the inequality in equation \eqref{eq:hormander-division} with
$\kappa=3N/2$ and $C_P=C a_0^{N/2}$.
\end{proof}

We next explain why Taylor truncation approximates the coefficients
in precisely the weighted norm supplied by the lemma. The essential
point is that its weight is invariant under a common rotation of all
coordinates. We use the Taylor truncations $\Pi_\ell$ defined in
Subsection~\ref{subsec:coefficient-norms}, applied componentwise to
vectors of functions.
Write $h=\sum_{n\ge0}h_n$ for the homogeneous Taylor expansion of an
entire function $h$.
For fixed $z$, Parseval's identity for
$h(e^{i\theta}z)=\sum_{n\ge0}e^{in\theta}h_n(z)$ gives
\[
 \frac1{2\pi}\int_0^{2\pi}|h(e^{i\theta}z)|^2\,d\theta
 =\sum_{n\ge0}|h_n(z)|^2.
\]
Integrate against $(1+|z|^2)^{-\kappa}d\gamma(z)$. This measure is unchanged
by $z\mapsto e^{i\theta}z$, so Tonelli's theorem identifies the squared
norm of $h$ with the sum of the squared norms of the $h_n$.
Whenever this norm is finite, applying the same identity to the tail gives
\begin{equation}\label{eq:weighted-taylor-convergence}
 \int |h-\Pi_\ell h|^2(1+|z|^2)^{-\kappa}\,d\gamma
 =\sum_{n>\ell}\int |h_n|^2(1+|z|^2)^{-\kappa}\,d\gamma
 \longrightarrow0.
\end{equation}

We can now prove uniqueness by induction on the number of variables.
For the base case, we set $R_0=\mathcal O(\C^0)=\C$, where
$\C^0$ consists of a single point.

\begin{proposition}\label{prop:affine-uniqueness}
 Let $I\subseteq R_d$ be an affine ideal. If
 $g\in\mathcal A(I)\cap\mathcal N(I)$ and $A_g(z) = e^{-|z|^2/2}|g(z)|$ has polynomial
 growth, then $g=0$.
\end{proposition}
\begin{proof}
We induct on $d$, allowing arbitrary polynomial growth exponent at each
step. The zero ideal gives $g=0$ by $\mathcal A(0)=\{0\}$, and the
unit ideal gives $g=0$ by $1^*(\partial)g=0$. This establishes the base cases, and establishes the statement for $d=0$. Suppose $d\ge1$
and $I$ is nonzero and proper.

First suppose that $A_g$ is rapidly decreasing. Then,
$\int |g|^2(1+|z|^2)^t\,d\gamma<\infty$ for every $t\ge0$.
Indeed, choose $k>t+d$, so the integrand relative to $d\mu$ is bounded by a
constant times $(1+|z|)^{2t-2k}$, which is integrable in real
dimension $2d$. In particular, $g\in\cF_d$.
Apply Lemma~\ref{lem:affine-limits} to the constant sequences
$H_\nu=H=I^h$ and $g_\nu=g$. It gives a representation $g=Pu$ with
$u$ a vector of entire functions, and $P=(p_1,\ldots,p_s)$ and each $p_i$ belongs to $I$.
Lemma~\ref{lem:gaussian-division} now gives
$g=\sum_i p_i h_i$ with
each $h_i$ square-integrable for $(1+|z|^2)^{-\kappa}d\gamma$.
Put $f_\ell=P[\Pi_\ell h]$ where $\Pi_{\ell} h$ is the $\ell$-th truncation of each entry of $h$. Each $f_{\ell}$ is a polynomial in $I$, so
$\langle g,f_\ell\rangle_{\cF_d}=0$ since
$g\in\mathcal N(I)\cap\cF_d=I^\perp$. Thus, since
$g-f_\ell=P(h-\Pi_\ell h)$, two applications of the Cauchy--Schwarz inequality yield
\[
 \begin{split}
 \|g\|_{\cF_d}^2 = \left|\int g\,\overline{g- f_{\ell}}\,d\gamma\right|
 &=\left|\int g\,\overline{P(h-\Pi_\ell h)}\,d\gamma\right|\\
 &\le
 \left(\int |g|^2|P|^2(1+|z|^2)^\kappa\,d\gamma\right)^{1/2}
 \left(\int |h-\Pi_\ell h|^2(1+|z|^2)^{-\kappa}\,d\gamma\right)^{1/2}
 \longrightarrow0.
 \end{split}
\]
The first factor is finite by rapid decrease and the polynomial growth of $P$. The second tends to zero by equation \eqref{eq:weighted-taylor-convergence}, applied to each component. Thus, $g=0$ in this case.

It remains to treat the case in which $A_g$ is not rapidly decreasing.
We will choose centers $b_\nu$ tending to infinity such that the
translates, normalized at those centers, have a common polynomial
bound on their weighted absolute values.
There is an integer $k\ge0$ for which $(1+|z|)^k A_g(z)$ is unbounded. Since $A_g$ has polynomial growth, there is some $m \geq 1$ such that $A_g(z)\le C(1+|z|)^m$. Fix an integer $q>\max\{k,m\}$. Put $\omega(z)=(1+|z|)^q$. The triangle inequality gives, for $x,y\in\C^d$,
\[
 \omega(x+y)\le\omega(x)\omega(y).
\]
Since $A_g\omega$ is unbounded, choose for each $\nu \in \mathbb{N}$ auxiliary points $a_\nu$ with
$A_g(a_\nu)\omega(a_\nu)\ge\nu$. These points locate large weighted
values of $A_g$. To obtain uniform bounds after translation, choose $b_\nu$
to maximize the continuous function
$x\mapsto A_g(x)/\omega(x-a_\nu)$. For each fixed $\nu$, this function
tends to zero as $|x|\to\infty$ because $q>m$, and its value at
$a_\nu$ is positive. Thus, its positive maximum is attained at some point $b_{\nu}$ and we also have $A_g(b_\nu)>0$.

The maximizing property has two consequences. First, comparing with
the value at $a=a_\nu$, and writing $b=b_\nu$, gives
\[
 \nu\le A_g(a)\omega(a)
 \le A_g(b)\frac{\omega(a)}{\omega(b-a)}
 \le A_g(b)\omega(b).
\]
The last inequality uses $a=b+(a-b)$ and $\omega(a-b)=\omega(b-a)$ together with the triangle inequality for $\omega$. Since $A_g\omega$ is unbounded, yet bounded on compact sets, it follows that
$|b_\nu|\to\infty$. Second, comparison with the value at $b+z$ gives
\[
 \frac{A_g(b+z)}{A_g(b)}
 \le\frac{\omega(b+z-a)}{\omega(b-a)}\le\omega(z).
\]
Consequently, the normalized weighted translates $g_\nu=\frac{D_{b_\nu}g}{(D_{b_\nu}g)(0)}$ are well defined, $g_\nu(0)=1$ and by equation \eqref{eq:weighted-translation} satisfy
\[
 A_{g_\nu}(z)=\frac{A_g(b_\nu+z)}{A_g(b_\nu)}\le\omega(z),
\]
In particular, they are uniformly bounded on compact sets. Montel's theorem gives a subsequential limit $G$ with $G(0)=1$ and $A_G\le\omega$.

Translation preservation puts
$g_\nu\in\mathcal A(I_{b_\nu})\cap\mathcal N(I_{b_\nu})$. Using Lemma~\ref{lem:translation-limits} to pass further to a subsequence, we get an affine ideal $L$ and a unit vector $e$ such that translation by $e$ maps $L$ into itself and $G\in\mathcal A(L)\cap\mathcal N(L)$. Up to applying a unitary change of coordinates, we may choose coordinates $(w,z')$ in which $e=(1,0,\ldots,0)$. Put $L'=L\cap\C[z']$ and $G'(z')=G(0,z')$. Lemma~\ref{lem:transverse-restriction} gives $G'\in\mathcal A(L')\cap\mathcal N(L')$. Restriction retains $G'(0)=1$, as well as the bound $e^{-|z'|^2/2}|G'(z')|\le(1+|z'|)^q$. The induction hypothesis in $d-1$ variables, with growth exponent $q$, gives $G'=0$. This contradicts $G'(0)=1$, completing the proof.
\end{proof}

\section{Singular vectors and the graded estimate}\label{sec:singular-vectors}

We now prove the uniform block bound required in the assumption of Proposition~\ref{prop:graded}. We first work with a scalar homogeneous
ideal. If its blocks were unbounded, singular vectors would give two
relations, one expressing orthogonality and the other ideal membership.
We will translate and normalize the vectors so that both relations
converge to the same entire function $G$, with $G(0)=1$ and $A_G\le1$.
Lemma~\ref{lem:translation-limits} will place $G$ in
$\mathcal A(L)\cap\mathcal N(L)$ for a limiting affine ideal $L$,
contradicting Proposition~\ref{prop:affine-uniqueness}.
We then pass from scalar ideals to modules of finite multiplicity by
encoding vector coordinates as auxiliary variables.

For a homogeneous ideal $I\subseteq R_d$, put
$I_n=I\cap R_{d,n}$ and $K_n=\cF_{d,n}\ominus I_n$. Recall the definition of $B_{j,n}$ from Section \ref{sec:preliminaries}, and fix, for a coordinate $j$, the abbreviation $B_n:=B_{j,n}$. Throughout the scalar argument the norms and block adjoints are Fock space norms and adjoints, given for $u\in K_n$ and $v\in I_{n+1}$ by 
\[
 B_nu=P_{I_{n+1}}(z_ju),\ \ \
 B_n^*v=P_{K_n}\partial_jv.
\]
A \emph{singular pair} for $B_n$ with singular value $\lambda>0$
consists of unit vectors $u\in K_n$ and $v\in I_{n+1}$ such that
$B_nu=\lambda v$ and $B_n^*v=\lambda u$.
When $B_n\ne0$, such a pair exists with $\lambda=\|B_n\|$:
take a unit eigenvector $u$ of $B_n^*B_n$ for its largest eigenvalue
$\lambda^2$, and put $v=B_nu/\lambda$.

Homogeneity also relates these polynomial spaces to the two spaces
of entire functions introduced in Subsection~\ref{subsec:affine-spaces}.
A polynomial in $K_m$ is orthogonal to $I_m$ by definition and to all
other degrees of $I$ by homogeneity. Thus,
\[
 K_m\subseteq\mathcal N(I),\qquad I_m\subseteq\mathcal A(I).
\]
These inclusions allow us to retain the two singular-vector relations
after translating the ideal.

\begin{proposition}\label{prop:scalar-block}
 For every homogeneous ideal $I\subseteq R_d$ and every $j$, we have $\sup_{n\ge0}\|B_n\|<\infty$.
\end{proposition}
\begin{proof}
 Suppose that $\lambda_\nu=\|B_{n_\nu}\|\to\infty$ for a sequence
 of degrees $n_\nu$. Since $\|B_n\|\le\sqrt{n+1}$, these degrees must
 tend to infinity. Temporarily suppress the index $\nu$ and write
 $n=n_\nu$ and $\lambda=\lambda_\nu>0$. Choose a singular pair
 $u\in K_n$, $v\in I_{n+1}$ with singular value $\lambda$.
 The projection formulas for $B_n$ and $B_n^*$ give
 \begin{equation}\label{eq:singular-residuals}
 z_ju-\lambda v\in K_{n+1},\qquad
 \partial_jv-\lambda u\in I_n.
 \end{equation}

 We first choose one translation and one normalization that control
 both $u$ and $v$. Let $a\in\C^d$ maximize $A_{u+v}$.
 This maximum is positive and is attained. Indeed, $u+v$ is a nonzero polynomial because its two nonzero summands have different degrees, and its weighted
 absolute value tends to zero at infinity. Their consecutive degrees
 also let us bound each summand by this maximum. Indeed,
 for $|\zeta|=1$,
 \[
 (u+v)(\zeta z)=\zeta^n(u(z)+\zeta v(z)).
 \]
 For each fixed $z$, a suitable $\zeta$ makes the triangle inequality
 an equality; if either summand vanishes, every $\zeta$ does so.
 Since $|\zeta z|=|z|$, this gives
 \[
 A_u(z)+A_v(z)
 =\max_{|\zeta|=1}A_{u+v}(\zeta z)\le A_{u+v}(a).
 \]
 Next, define 
 \[
 c=(D_a(u+v))(0),\qquad
 U=D_au/c,\qquad V=D_av/c,\qquad G=U+V.
 \]
 Then, $c\ne0$ and $|c|=A_{u+v}(a)$. The weighted translation
 identity \eqref{eq:weighted-translation} then yields for $z\in\C^d$, 
 \begin{equation}\label{eq:singular-normalization}
 G(0)=1,\qquad A_U(z)+A_V(z)\le1.
 \end{equation}

 We next use the residuals in equation \eqref{eq:singular-residuals} to approximate $G$ by a function in each
 of the two spaces associated with $I_a$. The weighted-translation
 product rules in equation \eqref{eq:translation-product-rules} give
 \begin{equation}\label{eq:singular-approximants}
 \begin{aligned}
 G-\frac{z_jU}{\lambda}
 &=\frac1{c\lambda}D_a\bigl((\lambda+a_j)u-(z_ju-\lambda v)\bigr)
 \in\mathcal N(I_a),\\
 G-\frac{\partial_jV}{\lambda}
 &=\frac1{c\lambda}D_a\bigl((\lambda+\overline{a_j})v
                         -(\partial_jv-\lambda u)\bigr)
 \in\mathcal A(I_a).
 \end{aligned}
 \end{equation}
 For the first membership, both $u$ and the first residual in equation
 \eqref{eq:singular-residuals} belong to $\mathcal N(I)$; for the
 second, both $v$ and the second residual belong to $\mathcal A(I)$.
 The assertions therefore follow from linearity and translation
 preservation. The terms $a_j U$ and $\overline{a_j}V$ introduced by the translation rules cancel against the corresponding scalar multiples. Hence, the remaining errors admit bounds on compact sets that are independent of the translating center $a$. Thus, before taking
 limits in $\nu$, the errors depend only on $U$, $V$, and $\lambda$.

 Restore now the index $\nu$. By equation \eqref{eq:singular-normalization},  the functions $G_\nu$ are uniformly bounded on every compact set.  Montel's theorem gives a locally uniform subsequential limit $G$  with $G(0)=1$ and $A_G\le1$. The two errors in equation \eqref{eq:singular-approximants} tend to zero in the same topology. More precisely, on $|z|\le R$ we have
 \[
 |z_jU_\nu(z)|\le R e^{R^2/2},\qquad
 |\partial_jV_\nu(z)|\le e^{(R+1)^2/2}.
 \]
 The second estimate is Cauchy's formula in the $j$th coordinate on
 a circle of radius one. Dividing these bounds by $\lambda_\nu\to\infty$ proves the required convergence of the error estimates to zero.

 Apply Lemma~\ref{lem:translation-limits} to $I$ and the sequence $(a_\nu)$. After passing to a further subsequence, we obtain a single affine ideal $L$ for both limits in equation \eqref{eq:singular-approximants}.  Since both functions on their left sides converge to $G$, we obtain  $G\in\mathcal N(L)\cap\mathcal A(L)$. Proposition~\ref{prop:affine-uniqueness} now gives $G=0$ because $A_G\le1$, contradicting $G(0)=1$.
\end{proof}

To obtain the finite-multiplicity result, we use Shalit's realization
of vector-valued polynomials as scalar polynomials linear in auxiliary
variables \cite[\S5, pp.~284--286]{shalitdivision}. The new variables record the vector coordinates; restricting to their degree-one part will recover the original module and its Fock blocks.
Let $t=(t_1,\ldots,t_r)$ be the auxiliary set of variables and let $e_1,\ldots,e_r$ be the standard
basis of $\C^r$. For a homogeneous $J\subseteq R_d^r$, define for $f=(f_1,\ldots,f_r)$, 
\[
 \Theta f=\sum_{a=1}^r t_af_a(z),\qquad
 I^\sharp=(\Theta f:f\in J)\subseteq\C[z,t].
\]
For $n,q\ge0$, let $E_{n,q}$ consist of polynomials homogeneous of
degree $n$ in $z$ and degree $q$ in $t$, equipped with the Fock norm
in $d+r$ variables. Such polynomials are \emph{bihomogeneous of
bidegree $(n,q)$}, and their total degree is $n+q$.
The total-degree-$k$ part of $I^\sharp$ is written $(I^\sharp)_k$. We now restate the main theorem and prove it using this realization.

\begin{maintheoremrepeat}
For every homogeneous submodule $J\subseteq R_d^r$ of infinite
codimension and every
$1\leq i,j\leq d$, there is a constant $C_{ij}$, independent of $n$,
such that for all $n\geq0$,
\[
\bigl\|[S_i^*,S_j]|_{\mathcal Q_{J,n}}\bigr\|
\leq\frac{C_{ij}}{n+1}.
\]
Writing $I=\Ann(R_d^r/J)$ and $\delta=\dim_\C V(I)$, the singular values satisfy $s_k([S_i^*,S_j])=O(k^{-1/\delta})$. Thus, we get that $[S_i^*,S_j]\in\cS_p$ for every $p>\delta$.
\end{maintheoremrepeat}
\begin{proof}
 By Proposition~\ref{prop:graded}, it suffices to prove that the
 Fock blocks $B_{j,n}$ for $J$ are bounded uniformly in $n$.
 We compare each of them with a scalar block for $I^\sharp$.
 Put $K_n=(\cF_{d,n}\otimes\C^r)\ominus J_n$.
 For each $n\ge0$, the map $\Theta$ is a Fock isometry from
 $\cF_{d,n}\otimes\C^r$ onto $E_{n,1}$. Indeed, the monomials
 $t_az^\alpha$ are pairwise orthogonal and have squared Fock norm
 $\alpha!$, which is also the squared norm of $z^\alpha e_a$.

 We first identify the submodule inside the scalar ideal. Since $J$ is homogeneous, $I^\sharp$ is generated by the
 $\Theta g$ with homogeneous $g\in J_m$, of bidegree $(m,1)$.
 Hence, $I^\sharp$ is the direct sum of its bidegree pieces,
 and we claim that
 \begin{equation}\label{eq:auxiliary-sector}
 (I^\sharp)_{n+1}\cap E_{n,1}=\Theta J_n.
 \end{equation}
 The inclusion from right to left follows from the definition.
 Conversely, write a member of the left side as a finite sum of
 polynomial multiples of the generators $\Theta g$, and take its
 $(n,1)$ component. When $g\in J_m$, only the multiplier's
 $(n-m,0)$ component contributes. That component is a polynomial
 in $z$ alone, and its product with $g$ lies in $J_n$. For $m>n$
 the component is zero. Applying $\Theta$ proves the other inclusion.

 Distinct bidegrees are Fock orthogonal. Thus, equation
 \eqref{eq:auxiliary-sector} identifies $\Theta K_n$ with the
 part of the scalar quotient graded component of total degree $n+1$ lying in $E_{n,1}$. In particular, $\Theta f\perp(I^\sharp)_{n+1}$ for
 $f\in K_n$. Likewise, $z_j\Theta f\in E_{n+1,1}$, so its
 projection onto $(I^\sharp)_{n+2}$ is its projection onto
 $\Theta J_{n+1}$. Since $\Theta$ is an isometry and commutes with
 multiplication by $z_j$, we obtain
 \[
 P_{(I^\sharp)_{n+2}}z_j\Theta f
 =\Theta\bigl(P_{J_{n+1}}z_jf\bigr)=\Theta B_{j,n}f.
 \]
 Proposition~\ref{prop:scalar-block}, applied to $I^\sharp$ in
 $d+r$ variables, bounds these scalar blocks independently of their
 total degree $n+1$. It follows that there exists a constant
 $C_{I^\sharp,j}>0$ independent of $n$ such that,
 \[
 \|B_{j,n}f\|_{\cF}
 =\bigl\|P_{(I^\sharp)_{n+2}}z_j\Theta f\bigr\|_{\cF}
 \le C_{I^\sharp,j}\|\Theta f\|_{\cF}
 =C_{I^\sharp,j}\|f\|_{\cF}.
 \]
 Therefore, $b_j<\infty$ for each original coordinate $z_j$.
 Proposition~\ref{prop:graded} now gives all the conclusions.
 Its singular-value count uses the original quotient grades $K_n$,
 whose growth is determined by
 $\delta=\dim_\C V(\Ann(R_d^r/J))$.
\end{proof}

We conclude by relating the quotient formulation to Arveson's original
question. The point is to pass to a tail whose defect space lies in
one degree, which makes the usual shift model homogeneous.

We say that commuting tuple $T$ is a \emph{row contraction} if
$\sum_jT_jT_j^*\le1$, and that it is \emph{pure} if the iterates of
$\Phi(A)=\sum_jT_jAT_j^*$ satisfy $\Phi^m(1)\to0$ strongly. We say it has
\emph{finite rank} if $D=1-\Phi(1)$ has finite rank. Theorem~\ref{thm:main} now implies that every graded pure commuting row
contraction of finite rank is $p$-essentially normal for $p>d$.

\begin{corollary} \label{cor:finite-rank-row-cont}
Let $T=(T_1,\ldots,T_d)$ be a pure, finite rank commuting row contraction on a Hilbert space $H$. Suppose that $H$ admits an orthogonal grading $H=\bigoplus_{n\in\mathbb Z}H_n$ such that $T_jH_n\subseteq H_{n+1}$ for $1\le j\le d,\ n\in\mathbb Z$. Then $[T_i^*,T_j]\in\cS_p$ for every $1\le i,j\le d$ and every $p>d$.
\end{corollary}

\begin{proof}
The finite-dimensional case is immediate, so suppose that $H$ is
infinite dimensional. By the construction in
\cite[proof of Theorem~B, pp.~220--221]{arveson2000}, there is an
integer $N$ such that $L=\bigoplus_{n\ge N}H_n$ has finite
codimension in $H$, the space $H_N$ is finite dimensional, and the
restriction $U=(T_1|_L,\ldots,T_d|_L)$ is pure with $1_L-\sum_{j=1}^d U_jU_j^*=P_{H_N}$.

Since $P_{H_N}v=v$ for $v\in H_N$, Arveson's dilation theorem
\cite[Theorem~1.4, p.~181]{arveson2000} gives a coisometric module
map $W:H_d^2\otimes H_N\to L$ satisfying for $\alpha\in\Nz^d,\ v\in H_N$,

\[
 W(z^\alpha\otimes v)=U^\alpha v\in H_{N+|\alpha|}.
\]

In particular, $W(M_{z_j}\otimes1)=U_jW$. Choose an orthonormal
basis of $H_N$, write $r=\dim H_N$, and let $J=\ker W\cap R_d^r$.
The intertwining identity makes $J$ a polynomial submodule. To see
that $J$ is homogeneous and that $\ker W=\overline J$, write
$f=\sum_{m\ge0}f_m\in\ker W$ as a sum of its homogeneous parts.
The monomial formula gives $Wf_m\in H_{N+m}$. These spaces are
mutually orthogonal, so $\sum_m Wf_m=0$ implies $Wf_m=0$ for every
$m$. Thus, every $f_m$ lies in $J$, and the polynomial partial sums
converge to $f$.

It follows that $W$ restricts to a unitary from $(H_d^2\otimes\C^r)\ominus\overline J$ onto $L$, intertwining the quotient shifts with $U$. Since $L$ is infinite dimensional and the images of the polynomials are dense in $L$, the algebraic quotient $R_d^r/J$ is infinite dimensional. Its annihilator variety is contained in $\C^d$ and hence has dimension at most $d$. Theorem~\ref{thm:main} therefore gives $[U_i^*,U_j]\in\cS_p$ for every $p>d$.

Finally, let $\widetilde U_j$ be the extension of $U_j$ by zero on
$L^\perp$. Then $\widetilde U_j=T_jP_L$, so $T_j-\widetilde U_j=T_j(1_H-P_L)$ has finite rank. It follows that $[T_i^*,T_j]-[\widetilde U_i^*,\widetilde U_j]$ also has finite rank. With respect to $H=L\oplus L^\perp$, we have $[\widetilde U_i^*,\widetilde U_j]=[U_i^*,U_j]\oplus 0_{L^\perp}.$ Since $[U_i^*,U_j]\in\cS_p$ for every $p>d$, the same then holds for $[T_i^*,T_j]$.
\end{proof}

\end{document}